\documentclass[12pt]{article}
\usepackage{amssymb}
\usepackage{amsmath}
\usepackage{array}
\usepackage[russian]{babel}
\usepackage[pdftex]{graphics}
\usepackage[X2,T2A]{fontenc}
\usepackage[cp1251]{inputenc}

\DeclareMathOperator{\RRe}{Re} \DeclareMathOperator{\IIm}{Im}
\DeclareMathOperator{\mmod}{mod}

\DeclareMathOperator{\vep}{\varepsilon}

\DeclareMathOperator{\vf}{\varphi}

\multlinegap=0em \DeclareFontFamily{T1}{msb}{}
\DeclareFontShape{T1}{msb}{m}{ol}{<5> <6> <7> <8> <9> gen * msbm
<10> <10.95> <12> <14.4> <17.28> <20.74> <24.88> msbm10}{}
\DeclareSymbolFont{AMSb}{T1}{msb}{m}{ol} \multlinegap=0em

\renewcommand{\S}{\mathhexbox278}
\renewcommand{\le}{\operatorname{\leqslant}}
\renewcommand{\ge}{\operatorname{\geqslant}}

\begin{document}

\begin{center}
{\rmfamily\bfseries\normalsize Kloosterman sums with primes and solvability of one congruence \\
with inverse residues - I}
\end{center}

\begin{center}
{\normalsize M.E.~Changa (Moscow), M.A.~Korolev\footnote{The second author was supported by the Russian Science Foundation under grant 19\,-11\,-00001.} (Moscow)}
\end{center}

\vspace{0.5cm}

\fontsize{11}{12pt}\selectfont

\textbf{Abstract.} In the paper, we establish a new estimate for Kloosterman sum over primes with respect to an arbitrary modulus $q$. This estimate together with some recent results of the second author are applied to the problem of solvability of the congruence
\[
g(p_{1})\,+\,\ldots\,+\,g(p_{k})\,\equiv\,m\pmod{q}
\]
in prime variables $p_{1},\ldots, p_{k}\leqslant N$, $N\le q^{1-\delta}$. Here $g(x)\,\equiv\,a\overline{x}+bx\pmod{q}$, where $a,b,m$ and $k\geqslant 3$ are arbitrary integers, $(ab,q)=1$.
The main result of the paper gives an asymptotic formula for the number of solutions for the case when $q$ is coprime to $6$.

\vspace{0.2cm}


\vspace{0.2cm}

Bibliography: 22 titles.

\fontsize{12}{15pt}\selectfont

\vspace{0.5cm}

\textbf{1. Introduction.}
\vspace{0.5cm}

Let $q\ge 2, a,b$ be integers. The exponential sum
\begin{equation}\label{1-01}
S(a,b;q)\,=\,\mathop{{\sum}'}\limits_{n=1}^{q}e_{q}(a\overline{n}+bn),
\end{equation}
is called a complete Kloosterman sum. Here $e_{q}(u) = e^{2\pi iu/q}$, the prime sign means the summation over $(n,q)=1$ and the symbol $\overline{n}=1/n$ stands for the inverse residue to $n$ modulo $q$, that is, for the solution of the congruence $n\overline{n} \equiv 1\pmod{q}$. The sum of the type (\ref{1-01}) where the summation is taken over some set which does not coincide with the reduced residual system modulo $q$, is called an incomplete Kloosterman sum. A particular case of such sum is Kloosterman sum with primes, that is
\begin{equation}\label{1-02}
W_{q}(a,b;X)\,=\,\mathop{{\sum}'}\limits_{p\le X}e_{q}(a\overline{p}+bp)
\end{equation}
Such sums were studied by E.~Fouvry, P.~Michel \cite{Fouvry_Michel_1998}, J.~Bourgain \cite{Bourgain_2005}, M.Z.~Garaev \cite{Garaev_2010}, E.~Fouvry, I.E.~Shparlinski \cite{Fouvry_Shparlinski_2011}, R.C.~Baker \cite{Baker_2012} and the author \cite{Korolev_2018a}-\cite{Korolev_2019a} (see also \cite{Bourgain_Garaev_2014}, \cite{Korolev_2017}).

An estimate of the sum (\ref{1-02}) for an arbitrary composite modulus $q$ is of special interest. In \cite{Fouvry_Shparlinski_2011}, \cite{Korolev_2018a}, such estimates for the case $b\equiv 0\pmod{q}$ (``homogeneous'' sum) are given. These estimates have some applications to various arithmetical problems (see, for example, \cite{Fouvry_Shparlinski_2011}, \cite{Baker_2012}, \cite{Korolev_2018c}).
In particular, such estimates allow one to study the solvability of the congruence
\begin{equation}\label{1-03}
\overline{p}_{1}\,+\,\ldots\,+\,\overline{p}_{k}\,\equiv\,m\pmod{q}
\end{equation}
in primes $p_{1},\ldots,p_{k}\le X$. In \cite{Korolev_2019b}, the second author find an estimate of ``inhomogeneous'' sum $W_{q}(a,b;X)$, $(ab,q)=1$, which is valid for $q^{\,3/4+\vep}\ll X\ll q^{\,3/2}$ (see the below Lemma 3.1). This estimate together with some other assertions allow one to study the solvability of the congruence
\begin{equation}\label{1-04}
g(p_{1})\,+\,\ldots\,+\,g(p_{k})\,\equiv\,m\pmod{q},\quad g(x)\,\equiv\,a\overline{x}+bx\pmod{q},
\end{equation}
in prime numbers $p_{1},\ldots,p_{k}\le X$. Such a problem is the main subject of the present paper.

The solution of this problem has some common features with circle method of Hardy, Littlewood and Ramanujan. Namely, the roles of ``major'' and ``minor'' arcs are played here by small and large divisors of
modulus $q$, respectively. First of them generate a main term, and the others contribute to the remainder. Consequently, an asymptotic formula for the number of solutions of (\ref{1-04}) involves a factor $\varkappa_{k}(a,b,m;q)$ which is similar to the singular series in the formulas given by circle method. In such interpretation, the estimate of the sum (\ref{1-02}) given in \cite{Korolev_2019b} corresponds to the estimate of the exponential sum over minor arcs. Consequently, in order to study (\ref{1-04}), we need some information about the behavior of (\ref{1-02}) over ``major arcs'', that is, for the case when the length $X$ of the sum is sufficiently large in comparison with modulus $q$. Such information is given by Theorems 2.1 and 2.2 below.

The main difficulty in the problem under considering is the study of the singular series $\varkappa_{k}(a,b,m;q)$ depending on the tuple of parameters $a, b, m$ and $k$. This singular series
is a multiplicative function of $q$ for fixed $a,b,m$ and $k$. It appears that if $k\ge 3$ is fixed and $q$ is coprime to $6$ then $\varkappa_{k}(a,b,m;q)$ is strictly positive for any tuple $(a,b,m)$ such that $1\le a,b,m\le q$, $(ab,q)=1$, and in this case the formula for the number $I_{k}(N)$ of solutions of (\ref{1-04}) given by Theorem 3.1 becomes asymptotic (see Theorem 7.1). At the same time, the behavior of $\varkappa_{k}(a,b,m;q)$ for $q = 2^{\,n}, 3^{\,n}$ is much more complicated. In particular, in these cases there exist tuples $(a,b,m)$ such that $\varkappa_{k}(a,b,m;q)=0$. That is the reason why we restrict here ourselves by the case $(q,6)=1$. The behavior of the values of $\varkappa_{k}(a,b,m;q)$, $q = 2^{\,n}, 3^{\,n}$ will be studied in a separate paper.

\textsc{Notations.} For integers $a$ and $b$, the symbol $(a,b)$ denotes their greatest common divisor, while the notation $(a;b)$ is used for the ordered pair. As usual, we use the standard notations $\varphi(n)$, $\Lambda(n)$, $\mu(n)$ and $\tau(n)$ for Euler, von Mangoldt, M\"{o}bius and the divisor functions, respectively. The number of different prime divisors of $n$ counted without multiplicity is denoted by $\omega(n)$.

\vspace{0.5cm}

\textbf{2. Estimates of ``long'' Kloosterman sums with primes.}

\vspace{0.5cm}

In this section, we prove several estimates for Kloosterman sum  $W_{q}(a,b;X)$ with primes in the case when the length $X$ of this sum is large in comparison with the modulus $q$.
Obviously, it is sufficient to estimate the sum
\[
T_{q}(X)\,=\,T_{q}(a,b;X)\,=\,\mathop{{\sum}'}\limits_{n\le X}\Lambda(n)e_{q}(a\overline{n}+bn).
\]
Below Theorems 2.1 and 2.2 are similar to Theorems 3 and 2 from \cite{Fouvry_Shparlinski_2011}.

We start with some auxiliary assertions.
\vspace{0.3cm}

\textsc{Lemma 2.1.}
\emph{Let $q\ge 2$, $a, b$ be integers. Then the following estimate holds:}
\[
\biggl|\,\mathop{{\sum}'}\limits_{n=1}^{q}e_{q}(a\overline{n}+bn)\biggr|\,\le\,\tau(q)\sqrt{q}\,(a,b,q)^{1/2}.
\]

\textsc{Corollary.} \emph{Under the assumptions of the lemma, for any $N\le q$ we have}
\[
\biggl|\,\mathop{{\sum}'}\limits_{n=1}^{N}e_{q}(a\overline{n}+bn)\biggr|\,\le\,\tau(q)\sqrt{q}\,(a,q)^{1/2}(\ln{q}+1).
\]

For the proof of Lemma 2.1, see \cite{Esterman_1961}. The corollary is derived from Lemma 2.1 by standard tools.
\vspace{0.2cm}

\textsc{Lemma 2.2.} \emph{Suppose that $X>2$ and let $a$ be even number, $2\le a\le X$. Then the number of primes $p\le X$ such that $p+a$ is also prime, does not exceed}
\[
\frac{cX}{(\ln{X})^{2\mathstrut}}\prod\limits_{p|a}\biggl(1+\frac{1}{p}\biggr),
\]
\emph{where $c>0$ is an absolute constant.}

\vspace{0.3cm}

For the proof, see, for example, \cite[Ch.II, \S 4, Th. 4.4]{Prachar_1957}.
\vspace{0.3cm}

\textsc{Lemma 2.3.} \emph{Suppose that $1<N<N_{1}\le 2N$, $1\le \delta\le N$. Then}
\[
\sum\limits_{\substack{N<n_{2}<n_{1}\le N_{1} \\ n_{1}-n_{2}\equiv 0\pmod{\delta}}}\Lambda(n_{1})\Lambda(n_{2})\,\ll\,N\sum\limits_{1\le m\le \tfrac{N}{\delta}}\;\sum\limits_{d|m\delta}\frac{\mu^{2}(d)}{d}\,+\,\frac{N^{\,3/2}}{\delta}\,\ln{N}.
\]

\textsc{Proof.} Denote by $W_{j}$ ($j = 1,2,3,4$) the contribution coming to the initial sum from the pairs $(n_{1};n_{2})$ of the form (a) $n_{1} = p_{1}$, $n_{2} = p_{2}$, (b) $n_{1} = p_{1}^{\ell}$, $n_{2} = p_{2}$, (c) $n_{1} = p_{1}$, $n_{2} = p_{2}^{r}$ and (d) $n_{1} = p_{1}^{\ell}$, $n_{2} = p_{2}^{r}$. Here $p_{1}, p_{2}$ are primes and $\ell,r$ are some integers with the conditions $2\le \ell,r\le H$, $H = [\log_{2}{N_{1}}]+1$.

If $p_{1} - p_{2} = m\delta$ then $1<m\delta\le N_{1}-N\le N$ and hence $1\le m\le N\delta^{-1}$. Thus, by Lemma 2.2, we get
\begin{multline*}
W_{1}\,\le\,(\ln{N_{1}})^{2}\sum\limits_{1\le m\le \tfrac{N}{\delta}}\;\sum\limits_{\substack{N<p_{2}\le N_{1} \\ p_{2}+m\delta\;\text{is prime}}}1\,\ll\,(\ln{N})^{2}
\sum\limits_{1\le m\le \tfrac{N}{\delta}}\frac{N}{(\ln{N})^{2}}\prod\limits_{p|m\delta}\biggl(1+\frac{1}{p}\biggr)\,\ll\\
\ll\,N\sum\limits_{1\le m\le \tfrac{N}{\delta}}\;\prod\limits_{p|m\delta}\biggl(1+\frac{1}{p}\biggr).
\end{multline*}

Next, estimating $W_{2}$, we fix an integer $\ell\ge 2$. Then
\begin{multline*}
\sum\limits_{\substack{N<p_{2}<p_{1}^{\ell}\le N_{1} \\ p_{1}^{\ell}-p_{2}\equiv 0\pmod{\delta}}}(\ln{p_{1}})(\ln{p_{2}})\,\ll\,(\ln{N})\sum\limits_{N^{1/\ell}<p_{1}\le N_{1}^{1/\ell}}(\ln{p_{1}})\sum\limits_{\substack{N<n\le p_{1}^{\ell} \\ n\equiv p_{1}^{\ell}\pmod{\delta}}}1\,\\
\ll\,\frac{N}{\delta}\,(\ln{N})\sum\limits_{N^{1/\ell}<p_{1}\le N_{1}^{1/\ell}}(\ln{p_{1}})\,\ll\,\frac{N^{1+1/\ell}}{\delta}.
\end{multline*}
The summation over $2\le \ell\le H$ yields:
\[
W_{2}\,\ll\,\frac{N^{3/2}}{\delta}\,\ln{N}.
\]
The same bound holds true for $W_{3}$. Finally, if we fix $\ell,r\ge 2$ then
\[
\sum\limits_{\substack{N<p_{2}^{r}<p_{1}^{\ell}\le N_{1} \\ p_{1}^{\ell}-p_{2}^{r}\equiv 0\pmod{\delta}}}(\ln{p_{1}})(\ln{p_{2}})\,\ll\,\sum\limits_{p_{1}\le N_{1}^{1/\ell}}\ln{p_{1}}
\sum\limits_{p_{2}\le N_{2}^{1/r}}\ln{p_{2}}\,\ll\,N^{1/\ell+1/r}.
\]
Summing over $2\le \ell,r\le H$ we find that $W_{4}\ll N$. Therefore,
\begin{multline*}
\sum\limits_{j=1}^{4}W_{j}\,\ll\,N\sum\limits_{1\le m\le \tfrac{N}{\delta}}\;\prod\limits_{p|m\delta}\biggl(1+\frac{1}{p}\biggr)\,+\,\frac{N^{3/2}}{\delta}\,\ln{N}+\,N\,\ll\\
\ll\,N\sum\limits_{1\le m\le \tfrac{N}{\delta}}\;\sum\limits_{d|m\delta}\frac{\mu^{2}(d)}{d}\,+\,\frac{N^{3/2}}{\delta}\,\ln{N}.\quad \Box
\end{multline*}
\vspace{0.3cm}

\textsc{Lemma 2.4.} \emph{Suppose that $1<M<M_{1}\le 2M$, $1<N<N_{1}\le 2N$, $MN\le X$, and let $|a_{m}|\le \ln{m}$, $|b_{m}|\le \tau(m)$ for $M<m\le M_{1}$. Then the sums}
\begin{align*}
& U(M,N)\,=\,\mathop{{\sum}'}\limits_{M<m\le M_{1}}a_{m}\mathop{{\sum}'}\limits_{\substack{N<n\le N_{1} \\ mn\le X}}e_{q}(a\overline{m}\overline{n}+bmn),\\
& W(M,N)\,=\,\mathop{{\sum}'}\limits_{M<m\le M_{1}}b_{m}\mathop{{\sum}'}\limits_{\substack{N<n\le N_{1} \\ mn\le X}}\Lambda(n)e_{q}(a\overline{m}\overline{n}+bmn)
\end{align*}
\emph{satisfy the estimates}
\begin{align}
& U(M,N)\,\ll\,\biggl(M\sqrt{N}+\frac{MN}{\sqrt[4\;]{q}}\,\tau(q)\,+\,N\sqrt{M}\sqrt[4\;]{q}\,\tau(q)\sqrt{\ln{X}}\biggr)\ln{X},\label{2-01}\\
& W(M,N)\,\ll\,\biggl(M\sqrt{N}+\frac{MN}{\sqrt[4\;]{q}}\,\frac{\tau(q)}{\sqrt{\ln{X}}}\,+\,N\sqrt{M}\sqrt[4\;]{q}\,\tau(q)\biggr)(\ln{X})^{2}. \label{2-02}
\end{align}

\textsc{Proof.} Writing both sums in a uniform manner as
\[
V(M,N)\,=\,\mathop{{\sum}'}\limits_{M<m\le M_{1}}\alpha(m)\mathop{{\sum}'}\limits_{\substack{N<n\le N_{1} \\ mn\le X}}\beta(n)e_{q}(a\overline{m}\overline{n}+bmn)
\]
and applying Cauchy's inequality we get
\begin{multline*}
|V(M,N)|^{2}\,\le\,\|\alpha\|_{2}^{2}\mathop{{\sum}'}\limits_{M<m\le M_{1}}\biggl|\;\mathop{{\sum}'}\limits_{\substack{N<n\le N_{1} \\ mn\le X}}\beta(n)e_{q}(a\overline{m}\overline{n}+bmn)\biggr|^{2},\\
\text{where}\quad \|\alpha\|_{2}\,=\,\biggl(\mathop{{\sum}'}\limits_{M<m\le M_{1}}|\alpha(m)|^{2}\biggr)^{\! 1/2}.
\end{multline*}
Further,
\begin{multline*}
|V(M,N)|^{2}\,\le\,\|\alpha\|_{2}^{2}\mathop{{\sum}'}\limits_{M<m\le M_{1}}\;\mathop{{\sum}'}\limits_{\substack{N<n_{1},n_{2}\le N_{1} \\ mn_{1}, mn_{2}\le X}}\beta(n_{1})\overline{\beta}(n_{2})e_{q}(a(\overline{n}_{1}-\overline{n}_{2})\overline{m}+b(n_{1}-n_{2})m)\,\le\\
\le\, \|\alpha\|_{2}^{2}\biggl\{\,\mathop{{\sum}'}\limits_{M<m\le M_{1}}\mathop{{\sum}'}\limits_{\substack{N<n\le N_{1} \\ mn\le X}}|\beta(n)|^{2}\,+\\
+\,2\mathop{{\sum}'}\limits_{N<n_{2}<n_{1}\le N_{1}}|\beta(n_{1})\beta(n_{2})|\,\biggl|\;\mathop{{\sum}'}\limits_{M<m\le M_{2}}e_{q}(a(\overline{n}_{1}-\overline{n}_{2})\overline{m}+b(n_{1}-n_{2})m)\biggr|\biggr\},
\end{multline*}
where $M_{2} = \min{\bigl(M_{1},Xn_{1}^{-1}\bigr)}$. Next, we split $(M,M_{2}]$ to the segments of length $q$ (with at most one possible exception, where this length is strictly less than $q$) and apply the estimates of Lemma 2.1 and its Corollary. Since $(\overline{n}_{1}-\overline{n}_{2},q) = (n_{1}-n_{2},q)$, we obtain
\[
|V(M,N)|^{2}\,\le\,M\|\alpha\|_{2}^{2}\|\beta\|_{2}^{2}\,+\,2\|\alpha\|_{2}^{2}\bigl(\bigl[Mq^{-1}\bigr]\sqrt{q}\tau(q)+\sqrt{q}\tau(q)(\ln{q}+1)\bigr)\,S,
\]
where
\[
S\,=\,\mathop{{\sum}'}\limits_{N<n_{2}<n_{1}\le N_{1}}|\beta(n_{1})\beta(n_{2})|\,(n_{1}-n_{2},q)^{1/2}\,\le\,\sum\limits_{\substack{\delta | q \\ \delta\le N}}\sqrt{\delta}
\mathop{{\sum}'}\limits_{\substack{N<n_{2}<n_{1}\le N_{1} \\ n_{1}-n_{2}\equiv 0\pmod{\delta}}}|\beta(n_{1})\beta(n_{2})|.
\]
For the case of $U(M,N)$, we obviously have
\[
\|\alpha\|_{2}^{2}\,\ll\,M(\ln{X})^{2},\quad \|\beta\|_{2}^{2}\,\ll\,N,\quad S\,\ll\,\sum\limits_{\substack{\delta | q \\ \delta\le N}}\sqrt{\delta}\,\frac{N^{2}}{\delta}\,\ll\,N^{2}\tau(q),
\]
and hence
\begin{multline*}
|U(M,N)|^{2}\,\ll\,M^{2}N(\ln{X})^{2}\,+\,M(\ln{X})^{2}\biggl(\frac{M}{\sqrt{q}}\,+\,\sqrt{q}\ln{q}\biggr)N^{2}\tau^{2}(q)\,\ll\\
\ll\,M^{2}N(\ln{X})^{2}\,+\,\frac{(MN)^{2}}{\sqrt{q}}\,(\ln{X})^{2}\tau^{2}(q)\,+\,MN^{2}\sqrt{q}(\ln{X})^{3}\tau^{2}(q),
\end{multline*}
so the desired estimate follows.

In the second case, we have
\[
\|\alpha\|_{2}^{2}\,\ll\,M(\ln{X})^{3},\quad \|\beta\|_{2}^{2}\,\ll\,N\ln{X}.
\]
Moreover, in view of Lemma 2.3 we have
\[
S\,\ll\,\sum\limits_{\substack{\delta|q \\ \delta\le N}}\sqrt{\delta}\biggl(N\sum\limits_{1\le m\le \tfrac{N}{\delta}}\;\sum\limits_{d|m\delta}\frac{\mu^{2}(d)}{d}\,+\,\frac{N^{3/2}}{\delta}\,\ln{N}\biggr)\,\ll\,N\Sigma\,+\,N^{3/2}(\ln{N})\tau(q),
\]
where
\[
\Sigma\,=\,\sum\limits_{\substack{\delta|q \\ \delta\le N}}\sqrt{\delta}\sum\limits_{1\le m\le \tfrac{N}{\delta}}\;\sum\limits_{d|m\delta}\frac{\mu^{2}(d)}{d}\,=\,
\sum\limits_{d\le N}\frac{\mu^{2}(d)}{d}\sum\limits_{\substack{\delta|q \\ \delta\le N}}\sqrt{\delta}\sum\limits_{\substack{1\le m\le \tfrac{N}{\delta} \\ m\delta\equiv 0\pmod{d}}}1.
\]
The inner sum does not exceed $N(d,\delta)(d\delta)^{-1}$. Hence,
\begin{multline*}
\Sigma\,\le\,\sum\limits_{d\le N}\frac{\mu^{2}(d)}{d}\sum\limits_{\substack{\delta|q \\ \delta\le N}}\sqrt{\delta}\,\frac{N}{d\delta}\,(d,\delta)\,=\\
=\,N\sum\limits_{d\le N}\frac{\mu^{2}(d)}{d^{2}}\sum\limits_{\substack{\delta|q \\ \delta\le N}}\frac{(d,\delta)}{\sqrt{\delta}}\,=\,N\sum\limits_{\substack{\delta|q \\ \delta\le N}}\frac{1}{\sqrt{\delta}}\sum\limits_{d = 1}^{+\infty}\frac{\mu^{2}(d)}{d^{2}}\,(d,\delta).
\end{multline*}
Obviously, the sum over $d$ is equal
\[
\prod\limits_{p}\biggl(1+\frac{(p,\delta)}{p^{2}}\biggr)\,=\,\prod\limits_{p\nmid{\delta}}\biggl(1+\frac{1}{p^{2}}\biggr)\prod\limits_{p|\delta}\biggl(1+\frac{1}{p}\biggr)\,\ll\,\prod\limits_{p|\delta}\biggl(1+\frac{1}{p}\biggr)
\,\ll\,\sum\limits_{d|\delta}\frac{\mu^{2}(d)}{d}.
\]
Thus,
\[
\Sigma\,\ll\,N\sum\limits_{\delta|q}\frac{1}{\sqrt{\delta}}\sum\limits_{d|\delta}\frac{\mu^{2}(d)}{d}\,\ll\,\sum\limits_{d|q}\frac{\mu^{2}(d)}{d}\sum\limits_{c|qd^{-1}}\frac{1}{\sqrt{cd}}\,\ll\,
\sum\limits_{d|q}\frac{\mu^{2}(d)}{d^{\,3/2}}\,\tau(qd^{-1})\,\ll\,\tau(q).
\]
Hence,
\[
S\,\ll\,N^{2}\tau(q)\,+\,N^{3/2}(\ln{N})\tau(q)\,\ll\,N^{2}\tau(q)
\]
and therefore
\begin{multline*}
|W(M,N)|^{2}\,\ll\,M^{2}N(\ln{X})^{4}\,+\,M(\ln{X})^{3}\biggl(\frac{M}{\sqrt{q}}+\sqrt{q}\ln{q}\biggr)N^{2}\tau^{2}(q)\,\ll\\
\ll\,M^{2}N(\ln{X})^{4}\,+\,\frac{(MN)^{2}}{\sqrt{q}}\,(\ln{X})^{3}\tau^{2}(q)\,+\,MN^{2}\sqrt{q}(\ln{X})^{4}\tau^{2}(q).
\end{multline*}
Now, the second inequality follows. $\Box$

\vspace{0.3cm}

\textsc{Theorem 2.1.}
\emph{Suppose that $0<\vep<0.1$ is arbitrary small but fixed constant, and let $q\ge q_{0}(\vep)$, $a, b$ be integers, $(ab,q)=1$. Then, for any $X\ge q\tau^{4}(q)(\ln{q})^{12}$, the following estimate holds}
\[
T_{q}(X)\ll X\Delta,\quad\textit{where}\quad \Delta\,=\,\frac{(\ln{X})^{5/2}}{\sqrt[4\;]{q}}\tau(q)\,+\,\Bigl(\frac{q}{X}\Bigr)^{\! 1/6}(\ln{X})^{2}(\tau(q))^{2/3}.
\]

\textsc{Proof.} Suppose that $1<V<\sqrt{X}$ and apply Vaughan's identity in the form given in \cite[Ch. II, \S 6, Th. 1]{Karatsuba_Voronin_1992}. Thus we get $T_{q}(X) = S_{1}-S_{2}-S_{3}-S_{4}+O(V)$, where
\begin{multline*}
S_{1}\,=\,\mathop{{\sum}'}\limits_{m\le V}\mu(m)\mathop{{\sum}'}\limits_{n\le X/m}(\ln{n})e_{q}(a\overline{m}\overline{n}+bmn),\quad S_{2}\,=\,\mathop{{\sum}'}\limits_{m\le V}a_{m}\mathop{{\sum}'}\limits_{n\le X/m}e_{q}(a\overline{m}\overline{n}+bmn),\\
S_{3}\,=\,\mathop{{\sum}'}\limits_{V<m\le V^{2}}a_{m}\mathop{{\sum}'}\limits_{n\le X/m}e_{q}(a\overline{m}\overline{n}+bmn),\quad S_{4}\,=\,\mathop{{\sum}'}\limits_{V<m\le X/V}b_{m}\mathop{{\sum}'}\limits_{V<n\le X/m}\Lambda(n)e_{q}(a\overline{m}\overline{n}+bmn),\\
a_{m}\,=\,\sum\limits_{r\ell = m,\;r,\ell\le V}\mu(r)\Lambda(\ell),\quad b_{m}\,=\,\sum\limits_{d|m,\, d\le V}\mu(d).
\end{multline*}
By Abel summation formula and Lemma 2.1, we write the inner sum $S_{1}(m)$ in $S_{1}$ as follows:
\begin{multline*}
S_{1}(m)\,=\,\mathop{{\sum}'}\limits_{n\le X/m}(\ln{n})e_{q}(a\overline{m}\overline{n}+bmn)\,=\,\mathbb{C}\bigl(X/m\bigr)\ln{\bigl(X/m\bigr)}\,-\,\int_{1}^{X/m}\mathbb{C}(u)\,\frac{du}{u},\\
\mathbb{C}(u)\,=\,\mathop{{\sum}'}\limits_{n\le u}e_{q}(a\overline{mn}+bmn)\,=\,\biggl[\frac{u}{q}\biggr]S(a,b;q)\,+\,\mathop{{\sum}'}\limits_{1\le n\le q\{u/q\}}e_{q}(a\overline{mn}+bmn)\,\ll\\
\ll\,\biggl[\frac{u}{q}\biggr]\sqrt{q}\tau(q)\,+\,\sqrt{q}\tau(q)\ln{q}\,\ll\,\biggl(\frac{u}{\sqrt{q}}\,+\,\sqrt{q}\ln{q}\biggr)\tau(q).
\end{multline*}
Hence,
\begin{multline*}
S_{1}(m)\,\ll\,\frac{X}{m}\,\frac{\tau(q)}{\sqrt{q}}\,\ln{X}\,+\,\sqrt{q}\tau(q)(\ln{q})\ln{X},\\
S_{1}\,\ll\,X\,\frac{\tau(q)}{\sqrt{q}}\,(\ln{X})^{2}\,+\,V\sqrt{q}\tau(q)(\ln{q})\ln{X}.
\end{multline*}
The same estimate holds for the sum $S_{2}$.

Next, we split the domains $(V,V^{\,2}]$ and $(1,XV^{-1}]$ of $m$ and $n$ in $S_{3}$ into segments of the type $M<m\le M_{1}$, $N<n\le N_{1}$ ($M_{1}\le 2M$, $N_{1}\le 2N$) by the points $M = 2^{s}V$, $N = XV2^{-t}$, $s,t = 1,2,\ldots$. Thus we obtain the expansion of $S_{3}$ into the sums $U(M,N)$ from Lemma 2.4. Using the inequality (\ref{2-01}) and denoting the summation over $M,N$ under considering by twin prime sign, we get
\begin{multline*}
S_{3}\,\ll\,(\ln{X})\mathop{{\sum}''}\limits_{V<M\le V^{2}}\mathop{{\sum}''}\limits_{1<N\le XM^{-1}}
\biggl(M\sqrt{N}+\frac{MN}{\sqrt[4\;]{q}}\,\tau(q)\,+\,N\sqrt{M}\sqrt[4\;]{q}\sqrt{\ln{X}}\biggr)\,\ll\\
\ll\, (\ln{X})\mathop{{\sum}''}\limits_{V<M\le V^{2}}\biggl(\sqrt{X}\sqrt{M}\,+\,\frac{X}{\sqrt[4\;]{q}}\,\tau(q)\,+\,\frac{X}{\sqrt{M}}\,\sqrt[4\;]{q}\sqrt{\ln{X}}\biggr)\,\ll\\
\ll\, V\sqrt{X}\ln{X}\,+\,\frac{X}{\sqrt[4\;]{q}}\,(\ln{X})^{2}\tau(q)\,+\,\frac{X}{\sqrt{V}}\,\sqrt[4\;]{q}(\ln{X})^{3/2}.
\end{multline*}
Similarly, splitting $(V,XV^{-1}]$ to diadic intervals and using (\ref{2-02}), we find that
\begin{multline*}
S_{4}\,\ll\,(\ln{X})^{2}\mathop{{\sum}''}\limits_{V<M\le XV^{-1}}\mathop{{\sum}''}\limits_{V<N\le XM^{-1}}
\biggl(M\sqrt{N}+\frac{MN}{\sqrt[4\;]{q}}\,\frac{\tau(q)}{\sqrt{\ln{X}}}\,+\,N\sqrt{M}\sqrt[4\;]{q}\tau(q)\biggr)\,\ll\\
\ll\, V\sqrt{X}(\ln{X})^{2}\,+\,\frac{X}{\sqrt[4\;]{q}}\,(\ln{X})^{5/2}\tau(q)\,+\,\frac{X}{\sqrt{V}}\,\sqrt[4\;]{q}(\ln{X})^{2}.
\end{multline*}
Summing the estimates for $S_{j}$, $1\le j\le 4$, we get
\[
T_{q}(X)\,\ll\, V\sqrt{q}(\ln{X})^{2}\tau(q)\,+\,V\sqrt{X}(\ln{X})^{2}\,+\,\frac{X}{\sqrt[4\;]{q}}\,(\ln{X})^{5/2}\tau(q)\,+\,\frac{X}{\sqrt{V}}\,\sqrt[4\;]{q}(\ln{X})^{2}\tau(q).
\]
Now we choose $V$ balancing the second and the last terms:
\[
V\sqrt{X}(\ln{X})^{2}\,=\,\frac{X}{\sqrt{V}}\,\sqrt[4\;]{q}(\ln{X})^{2}\tau(q),\quad \quad\text{that is,}\quad V\,=\,X^{1/3}q^{1/6}(\tau(q))^{2/3}.
\]
The condition $V\le \sqrt{X}$ is equivalent to the inequality $X\ge q\tau^{4}(q)$ which is obviously true for $X$ under considering. Thus we obtain:
\begin{multline*}
T_{q}(X)\,\ll\,X^{1/3}q^{2/3}(\ln{X})^{2}(\tau(q))^{5/3}\,+\,X^{5/6}q^{1/6}(\ln{X})^{2}(\tau(q))^{2/3}\,+\,\frac{X}{\sqrt[4\;]{q}}\,(\ln{X})^{5/2}\tau(q)\,\ll\\
\ll\,X\biggl\{\biggl(\frac{q}{X}\,(\tau(q))^{5/2}(\ln{X})^{3}\biggr)^{\! 2/3}\,+\,\biggl(\frac{q}{X}\,(\tau(q))^{4}(\ln{X})^{12}\biggr)^{\! 1/6}\,+\,\frac{(\ln{X})^{5/2}}{\sqrt[4\;]{q}}\,\tau(q)\biggr\}\,\ll\\
\ll\,X\biggl\{\biggl(\frac{q}{X}\biggr)^{\! 1/6}(\tau(q))^{2/3}(\ln{X})^{2}\,+\,\frac{(\ln{X})^{5/2}}{\sqrt[4\;]{q}}\,\tau(q)\biggr\}.
\end{multline*}
Theorem is proved. $\Box$
\vspace{0.3cm}

\textsc{Corollary.} \emph{Under the conditions of Theorem} 2.1, $W_{q}(a,b;X)\ll \pi(X)\Delta$.
\vspace{0.3cm}

\textsc{Proof.} By Abel summation formula,
\[
W_{q}(a,b;X)\,=\,\mathop{{\sum}'}\limits_{n\le X}\frac{\Lambda(n)}{\ln{n}}\,e_{q}(a\overline{n}+bn)\,+\,O(\sqrt{X})\,=\,\frac{T_{q}(X)}{\ln{X}}\,+\,\int_{2}^{X}\frac{T_{q}(u)du}{u(\ln{u})^{2}}\,+\,O(\sqrt{X}).
\]
Setting $Y = q\tau^{4}(q)(\ln{q})^{12}$, we estimate $T_{q}(u)$ trivially for $2\le u\le Y$ and by Theorem 2.1 for $Y<u\le X$. Thus we get
\begin{multline*}
W_{q}(a,b;X)\,\ll\,\pi(X)\Delta\,+\,\frac{Y}{(\ln{Y})^{2}}\,+\,\int_{Y}^{X}\Bigl(u^{-1/6}q^{1/6}(\tau(q))^{2/3}+\frac{\tau(q)}{\sqrt[4\;]{q}}\sqrt{\ln{u}}\Bigr)du
\ll\\
\ll\,\pi(X)\Delta\,+\,\frac{Y}{(\ln{Y})^{2}}\,+\,\frac{X\Delta}{(\ln{X})^{2}}\,\ll\,\pi(X)\Delta.
\end{multline*}
Corollary is proved. $\Box$

The Generalized Riemann hypothesis allows one to improve the estimates of Theorem 2.1 and its corollary. Namely, the following assertion holds true.
\vspace{0.3cm}

\textsc{Theorem 2.2.} \emph{Let $q, a, b$ be integers, $(ab,q)=1$ and suppose that $X\ge q(\ln{q})^{4}$. Then, under the generalized Riemann hypothesis, the following estimate holds:}
\[
T_{q}(X)\ll X\Delta_{1},\quad\textit{where}\quad \Delta_{1}\,=\,\frac{\sqrt{q}}{\varphi(q)}\tau(q)\,+\,\sqrt{\frac{q}{X}}\,(\ln{X})^{2}.
\]
\textsc{Proof.} First we write $T_{q}(X)$ as follows:
\[
T_{q}(X)\,=\,\mathop{{\sum}'}\limits_{\nu = 1}^{q}\sum\limits_{\substack{n\le X \\ n\equiv\nu \pmod{q}}}\Lambda(n)e_{q}(a\overline{n}+bn)\,=\,
\mathop{{\sum}'}\limits_{\nu = 1}^{q}e_{q}(a\overline{\nu}+b\nu)\sum\limits_{\substack{n\le X \\ n\equiv\nu \pmod{q}}}\Lambda(n).
\]
Further, the inner sum is written as
\[
\sum\limits_{n\le X}\frac{\Lambda(n)}{\varphi(q)}\sum\limits_{\chi\pmod{q}}\overline{\chi}(n)\chi(\nu)\,=\,\frac{1}{\varphi(q)}\sum\limits_{\chi\mmod{q}}\chi(\nu)\sum\limits_{n\le X}
\overline{\chi}(n)\Lambda(n)\,=\,\frac{1}{\varphi(q)}\sum\limits_{\chi\mmod{q}}\chi(\nu)\psi(X;\overline{\chi}),
\]
where
\[
\psi(X;\chi)\,=\,\sum\limits_{n\le X}\chi(n)\Lambda(n).
\]
Thus,
\[
T_{q}(X)\,=\,\frac{1}{\varphi(q)}\sum\limits_{\chi\mmod{q}}\psi(X;\overline{\chi})\mathop{{\sum}'}\limits_{\nu = 1}^{q}\chi(\nu)e_{q}(a\overline{\nu}+b\nu)\,=\,
\frac{1}{\varphi(q)}\sum\limits_{\chi\mmod{q}}\psi(X;\overline{\chi})S(a,b,q;\chi),
\]
where
\[
S(a,b,q;\chi)\,=\,\mathop{{\sum}'}\limits_{\nu = 1}^{q}\chi(\nu)e_{q}(a\overline{\nu}+b\nu).
\]
Extracting  the contribution from the principal character $\chi_{0}$ we find that
\[
T_{q}(X)\,=\,\frac{\psi(X;\chi_{0})}{\vf(q)}\,S(a,b;q)\,+\,R,\quad\;\;\text{where}\;\;
|R|\,\le\,\frac{1}{\varphi(q)}\sum\limits_{\chi\ne\chi_{0}}|\psi(X;\chi)|\cdot|S(a,b,q;\chi)|.
\]
One can check that
\[
\psi(X;\chi_{0})\,=\,\psi(X)\,+\,O((\ln{X})^{2})\,\ll\,X.
\]
Moreover, the Generalized Riemann hypothesis yields that $\psi(X;\chi)\ll \sqrt{X}(\ln{X})^{2}$ for any $\chi\ne\chi_{0}$ (see, for example, \cite[Ch. 20]{Davenport_1980}). Hence,
\[
T_{q}(X)\,\ll\,\frac{X}{\varphi(q)}\,\sqrt{q}\tau(q)\,+\,\frac{\sqrt{X}}{\varphi(q)}\,(\ln{X})^{2}S,
\]
where
\[
S\,=\,\sum\limits_{\chi\mmod{q}}|S(a,b;q,\chi)|.
\]
By Cauchy's inequality,
\begin{multline*}
S^{2}\,\le\,\varphi(q)\sum\limits_{\chi\mmod{q}}\biggl|\mathop{{\sum}'}\limits_{\nu = 1}^{q}\chi(\nu)e_{q}(a\overline{\nu}+b\nu)\biggr|^{2}\,=\,\varphi(q)\sum\limits_{\chi\mmod{q}}\mathop{{\sum}'}\limits_{\mu,\nu = 1}^{q}\chi(\nu)\overline{\chi}(\mu)e_{q}(a(\overline{\nu}-\overline{\mu})+b(\nu-\mu))\,=\\
=\,\varphi(q)\mathop{{\sum}'}\limits_{\mu,\nu = 1}^{q}e_{q}(a(\overline{\nu}-\overline{\mu})+b(\nu-\mu))
\sum\limits_{\chi\mmod{q}}\chi(\nu)\overline{\chi}(\mu)\,=\,\varphi^{3}(q).
\end{multline*}
Therefore,
\[
T_{q}(X)\,\ll\,\frac{X}{\varphi(q)}\,\sqrt{q}\,\tau(q)\,+\,\sqrt{X}\sqrt{\varphi(q)}\,(\ln{X})^{2}\,\ll\,X\Delta_{1}. \quad \Box
\]
\textsc{Corollary.} \emph{Under the condition of Theorem} 2.2, $W_{q}(a,b;X)\ll \pi(X)\Delta_{1}$.
\vspace{0.3cm}

The proof follows the same lines as the proof of Corollary to Theorem 2.1, but with $Y = q(\ln{q})^{4}$.
\vspace{0.5cm}

\textbf{3. Application to a congruence with inverses.}

\vspace{0.5cm}

Following two lemmas are proved in \cite{Korolev_2019b}.
\vspace{0.3cm}

\textsc{Lemma 3.1.} \emph{Let $0<\vep<0.1$ be any fixed constant, $q\ge q_{0}(\vep)$, $(ab,q)=1$. Then, for any $X$ satisfying the conditions $q^{\,3/4+\vep}\le X\le\,(q/2)^{3/2}$, the sum
\[
W_{q}(X)\,=\,W_{q}(a,b;X)\,=\,\mathop{{\sum}'}\limits_{p\le X}e_{q}(a\overline{p}+bp)
\]
obeys the estimate $W_{q}(X)\ll Xq^{\vep}\Delta$, where}
\begin{equation*}
\Delta\,=\,
\begin{cases}
\bigl(q^{\,3/4}X^{-1}\bigr)^{1/7}, & \textit{when}\quad q^{\,3/4}\le X\le q^{\,7/8},\\
\bigl(q^{\,2/3}X^{-1}\bigr)^{3/35}, & \textit{when}\quad q^{\,7/8}\le X\le\,(q/2)^{3/2}.
\end{cases}
\end{equation*}
\vspace{0.3cm}

\textsc{Lemma 3.2.} \emph{Suppose that $(ab,q)=1$. Then the number of solutions of the congruence $g(x)\equiv g(y)\pmod{q}$ with the conditions $1\le x,y\le q$ does not exceed} $\kappa(q)\le 2^{\omega(q)+1}\tau(q)q$.
\vspace{0.1cm}

These assertions together with the estimates of Theorems 2.1, 2.2 allow one to study the solvability of some congruences to any composite modulus. In what follows, we prove
\vspace{0.2cm}

\textsc{Theorem 3.1.} \emph{Let $0<\vep<0.01$ be an arbitrary fixed constant and let $k\ge 3$ be any fixed integer. Suppose that $q\ge q_{0}(\vep, k)$. Further, let $(ab,q)=1$ and $g(x)\equiv a\overline{x}+bx \pmod{q}$. Finally, let}
\[
c_{k}\,=\,\frac{2(k+33)}{3k+64}\quad \textit{if}\quad 3\le k\le 16\quad\textit{and}\quad c_{k}\,=\,\frac{3k+50}{4(k+12)}
\quad \textit{if}\quad k\ge 17,
\]
\emph{and suppose that $q^{\,c_{k}+\vep}\le N\le q$. Then the number $I_{k}(N) = I_{k}(N,q,a,b,m)$ of solutions of} (\ref{1-04}) \emph{in primes $p_{j}\le N$, $(p_{j},q)=1$, is equal to}
\[
I_{k}(N)\,=\,\frac{\pi^{k}(N)}{q}\,\bigl(\varkappa_{k}(q)\,+\,O(\Delta_{k})\bigr).
\]
\emph{Here $\varkappa_{k}(q) = \varkappa_{k}(a,b,m;q)$ is some non-negative multiplicative function of $q$ for any fixed tuple $k,a,b$ and $m$. Moreover,}

a) \emph{for any $k\ge 7$ we have $\Delta_{k} = (\ln\ln{N})^{B}(\ln{N})^{-A}$, where}
\[
A\,=\,\frac{1}{2}+\frac{29}{2}(k-7),\quad B\,=\,2^{k}-1;
\]

b) \emph{for any $k\ge 3$ we have $\Delta_{k} = q^{\,-\vep}$, if Generalized Riemann hypothesis is true.}
\vspace{0.1cm}

\textsc{Remark.} For any $k\ge 3$, we have $c_{k}\le 1-\tfrac{1}{73}<1$.
\vspace{0.1cm}

\textsc{Proof.} Obviously, we have
\[
I_{k}(N)\,=\,\frac{1}{q}\sum\limits_{c = 1}^{q}e_{q}(-cm)W_{q}^{k}(ac,bc;N).
\]
If $(c,q)=\delta$ then $q = \delta r$, $c = \delta f$ for some $f,r$ with $(r,f)=1$. Hence,
\begin{equation}\label{3-01}
I_{k}(N)\,=\,\frac{1}{q}\sum\limits_{r|q}\sum\limits_{\substack{f = 1 \\ (f,r)=1}}^{r}e_{r}(-mf)W_{\delta r}^{k}(a\delta f,b\delta f;N).
\end{equation}
For the proof of the unconditional result, we set $\mathcal{L} = \ln{N}$, $F = \mathcal{L}^{C}$, $C = 70$, $G = 2N^{2/3}$ and split the sum in (\ref{3-01}) into three parts $I_{k}^{(1)}$, $I_{k}^{(2)}$ and $I_{k}^{(3)}$ satisfying the intervals $1\le r\le F$, $F<r\le G$, $G<r\le q$ (some of them may be empty). The first part contributes to the main term in the asymptotic for $I_{k}(N)$, and the others give the remainder.

\textsc{The calculation of} $I_{k}^{(1)}$. First, we note that
\begin{equation}\label{3-02}
W_{\delta r}(a\delta f,b\delta f;N)\,=\,\sum\limits_{\substack{p\le N \\ p\nmid q}}e_{r}(f(a\overline{p}+bp))\,=\,W_{r}(af,bf;N)\,+\,\theta_{1}\omega(q),\quad |\theta_{1}|\le 1.
\end{equation}
Further,
\[
W_{r}(af,bf;N)\,=\,\sum\limits_{\substack{h = 1 \\ (h,r)=1}}^{r}e_{r}(f(a\overline{h}+bh))\sum\limits_{\substack{p\le N \\ p\equiv h\pmod{r}}}1\,=\,
\sum\limits_{\substack{h = 1 \\ (h,r)=1}}^{r}e_{r}(f(a\overline{h}+bh))\pi(N;r,h).
\]
By Siegel-Walfisz theorem, for $1<r\le F$ we have
\[
\pi(N;r,h)\,=\,\frac{\pi(N)}{\vf(r)}\,+\,O\bigl(\pi(N)e^{-c\sqrt{\mathcal{L}}}\bigr),\quad c>0.
\]
Hence,
\[
W_{\delta r}(a\delta f,b\delta f;N)\,=\,\frac{\pi(N)}{\varphi(r)}\,S(af,bf;r)\,+\,O\bigl(r\pi(N)e^{-c\sqrt{\mathcal{L}}}\bigr).
\]
Consequently,
\begin{multline*}
W_{\delta r}^{k}(a\delta f,b\delta f;N)\,=\\
=\,\frac{\pi^{k}(N)}{\varphi^{k}(r)}\,S^{k}(af,bf;r)\,+\,
O\biggl(\,\frac{\pi^{k}(N)}{\varphi^{k-1}(r)}\,|S(af,bf;r)|^{k-1}e^{-c\sqrt{\mathcal{L}}}\biggr)\,+\,
O\bigl(\pi^{k}(N)r^{k}e^{-ck\sqrt{\mathcal{L}}}\bigr)\,=\\
=\,\frac{\pi^{k}(N)}{\varphi^{k}(r)}\,S^{k}(af,bf;r)\,+\,O\bigl(\pi^{k}(N)e^{-c_{1}\sqrt{\mathcal{L}}}\bigr),
\quad c_{1} = \tfrac{2}{3}\,c.
\end{multline*}
Hence, turning to the sum $I_{k}^{(1)}$, we get
\begin{multline*}
I_{k}^{(1)}\,=\,\frac{1}{q}\sum\limits_{\substack{r|q \\ 1\le r\le F}}e_{r}(-mf)\biggl\{\frac{\pi^{k}(N)}{\varphi^{k}(r)}\,S^{k}(af,bf;r)\,+\,
O\bigl(\pi^{k}(N)e^{-c_{1}\sqrt{\mathcal{L}}}\bigr)\biggr\}\,=\\
=\,\frac{\pi^{k}(N)}{q}\sum\limits_{\substack{r|q \\ 1\le r\le F}}A_{k}(r)\,+\,O\biggl(\,\frac{\pi^{k}(N)}{q}\,e^{-c_{2}\sqrt{\mathcal{L}}}\bigg),
\end{multline*}
where
\[
A_{k}(r)\,=\,A_{k}(a,b,m;r)\,=\,\frac{1}{\varphi^{k}(r)}\sum\limits_{\substack{f=1 \\ (f,r)=1}}^{r}e_{r}(-mf)\,S^{k}(af,bf;r),\quad c_{2} = \tfrac{1}{3}\,c.
\]
Now we show that $A_{k}(r)$ is a multiplicative function of $r$ (for any fixed $a,b,m$ and $k$). Indeed,
let $(r_{1},r_{2})=1$, where $(r_{1},r_{2})=1$ and $r_{1}, r_{2}>1$. If $y$ and $z$ run through the reduced residual systems moduli $r_{1}$ and $r_{2}$ then $x = yr_{2}+zr_{1}$ runs through the reduced residual system modulo $r$ and, moreover,
\[
\overline{x}\,\equiv\,\overline{y}\,\overline{r}_{2}^{2}r_{2}\,+\,\overline{z}\,\overline{r}_{1}^{2}r_{1}\pmod{r},
\]
where $y\overline{y}\equiv r_{2}\overline{r}_{2}\equiv 1\pmod{r_{1}}$, $z\overline{z}\equiv r_{1}\overline{r}_{1}\equiv 1 \pmod{r_{2}}$. Therefore,
\begin{multline*}
S(u,v;r)\,=\,\sum\limits_{\substack{x=1 \\ (x,r)=1}}^{r}\exp{\biggl(\frac{2\pi i}{r}\,(u\overline{x}+vx)\biggr)}\,=\,\\
=\,\sum\limits_{\substack{y=1 \\ (y,r_{1})=1}}^{r_{1}}
\sum\limits_{\substack{z=1 \\ (z,r_{2})=1}}^{r_{2}}\exp{\biggl(\frac{2\pi i}{r_{1}r_{2}}\,(u\overline{y}\,\overline{r}_{2}^{2}r_{2}+u\overline{z}\,\overline{r}_{1}^{2}r_{1}+vyr_{2}+vzr_{1})\biggr)}\,=\\
=\,\sum\limits_{\substack{y=1 \\ (y,r_{1})=1}}^{r_{1}}\exp{\biggl(\frac{2\pi i}{r_{1}}\,(u\overline{y}\,\overline{r}_{2}^{2}+vy)\biggr)}\,
\sum\limits_{\substack{z=1 \\ (z,r_{2})=1}}^{r_{2}}\exp{\biggl(\frac{2\pi i}{r_{2}}\,(u\overline{z}\,\overline{r}_{1}^{2}+vz)\biggr)}.
\end{multline*}
Changing the variables $y, z$ to $y\overline{r}_{2}$, $z\overline{r}_{1}$, we get
\begin{multline*}
S(u,v;r)\,=\,\sum\limits_{\substack{y=1 \\ (y,r_{1})=1}}^{r_{1}}\exp{\biggl(\frac{2\pi i\overline{r}_{2}}{r_{1}}\,(u\overline{y}+vy)\biggr)}
\sum\limits_{\substack{z=1 \\ (z,r_{2})=1}}^{r_{2}}\exp{\biggl(\frac{2\pi i\overline{r}_{1}}{r_{2}}\,(u\overline{z}\,+vz)\biggr)}\,=\\
=\,S(u\overline{r}_{2},v\overline{r}_{2};r_{1})
S(u\overline{r}_{1},v\overline{r}_{1};r_{2}).
\end{multline*}
Further, if $f$ runs through the reduced residual system modulo $r = r_{1}r_{2}$ then $f$ is represented in the form $f = sr_{2}+tr_{1}$ where $s,t$ run through the reduced systems moduli $r_{1}$ and $r_{2}$. Thus, we obtain
\begin{multline*}
S(af,bf;r)\,=\,S(af\overline{r}_{2},bf\overline{r}_{2};r_{1})\,S(af\overline{r}_{1},bf\overline{r}_{1};r_{2})\,=\,
S(as,bs;r_{1})S(at,bt;r_{2}),\\
e_{r}(-mf)\,=\,e_{r_{1}}(-ms)e_{r_{2}}(-mt).
\end{multline*}
Therefore,
\begin{multline*}
A_{k}(r)\,=\,\frac{1}{\varphi^{k}(r)}\sum\limits_{\substack{f=1 \\ (f,r)=1}}^{r}e_{r}(-mf)S^{k}(af,bf;r)\,=\\
=\,\frac{1}{\varphi^{k}(r_{1})}\sum\limits_{\substack{s=1 \\ (s,r_{1})=1}}^{r_{1}}e_{r_{1}}(-ms)S^{k}(as,bs;r_{1})\cdot \frac{1}{\varphi^{k}(r_{2})}\sum\limits_{\substack{t=1 \\ (t,r_{2})=1}}^{r_{2}}e_{r_{2}}(-mt)S^{k}(at,bt;r_{2})\,=\\
=\,A_{k}(r_{1})A_{k}(r_{2}).
\end{multline*}
If $(u,v,r)=1$ then, by Lemma 2.1,
\begin{equation}\label{3-03}
|A_{k}(r)|\,<\,\frac{\varphi(r)\,r^{\,k/2}\tau^{k}(r)}{\varphi^{k}(r)}\,\ll_{\vep}\,r^{\,-k/2+1+\vep},
\end{equation}
and hence
\[
\biggl|\sum\limits_{\substack{F<r\le q \\ r|q}}A_{k}(r)\biggr|\,\ll_{\vep}\,\sum\limits_{r>F}r^{\,-k/2+1+\vep}\,\ll_{\vep}\,F^{\,-k/2+2+\vep}
\]
for $k\ge 2$. Therefore,
\begin{multline}\label{3-04}
I_{k}^{(1)}\,=\,\frac{\pi^{k}(N)}{q}\,\varkappa_{k}(q)\,+\,O_{\vep}\biggl(\frac{\pi^{k}(N)}{q}
\bigl(F^{\,-k/2+2+\vep}\,+\,e^{-c_{2}\sqrt{\mathcal{L}}}\bigr)\biggr)\,=\\
=\,\frac{\pi^{k}(N)}{q}\,\varkappa_{k}(q)\,+\,O_{\vep}\biggl(\frac{\pi^{k}(N)}{q}\,\mathcal{L}^{-A(k/2-2-\vep)}\biggr),
\end{multline}
where
\[
\varkappa_{k}(q)\,=\,\varkappa_{k}(a,b,m;q)\,=\,\sum\limits_{r|q}A_{k}(q).
\]
Obviously, $\varkappa_{k}(q)$ is multiplicative function for fixed $a,b,m$ and $k$.

\textsc{The estimation of} $I_{k}^{(j)}$, $j = 2,3$. Denote by $E_{j}$, $j = 2,3$ the segments $F<r\le G$ and $G<r\le q$. Then, by (\ref{3-02}) we have
\begin{multline*}
|I_{k}^{(j)}|\,\le\,\frac{1}{q}\sum\limits_{\substack{r|q \\ r\in E_{j}}}\sum\limits_{\substack{f=1 \\ (f,r)=1}}^{r}\bigl|W_{r}(af,bf;N)\,+\,\theta_{1}\omega(q)\bigr|^{k}\,\le\\
\le\,\frac{2^{k-1}}{q}\sum\limits_{\substack{r|q \\ r\in E_{j}}}\sum\limits_{\substack{f=1 \\ (f,r)=1}}^{r}\bigl(|W_{r}(af,bf;N)|^{k}\,+\,\omega^{k}(q)\bigr)\,\le\\
\le\,\frac{2^{k-1}}{q}\sum\limits_{\substack{r|q \\ r\in E_{j}}}W_{r}^{k-2}\sum\limits_{f=1}^{r}|W_{r}(af,bf;N)|^{2}\,+\,\frac{(2\omega(q))^{k}}{q}\,\sum\limits_{r|q}\varphi(r)\,=\\
=\,\frac{2^{k-1}}{q}\sum\limits_{\substack{r|q \\ r\in E_{j}}}W_{r}^{k-2}\,rJ_{r}(N)\,+\,(2\omega(q))^{k},
\end{multline*}
where $\displaystyle W_{r}=\max_{(f,r)=1}|W_{r}(af,bf;N)|$ and $J_{r}(N)$ denotes the number of solutions of the congruence
\[
g(p_{1})\,\equiv\,g(p_{2})\pmod{r}
\]
in primes $1<p_{1},p_{2}\le N$. Splitting the domain of $p_{1}$ and $p_{2}$ into segments of length $r$ and using Lemma 3.2 we get:
\[
rJ_{r}(N)\,\le\,\bigl(\bigl[Nr^{-1}\bigr]+1\bigr)^{2}2^{\,\omega(r)+1}r^{2}\tau(r)\,\ll\,(N^{2}+r^{2})\tau^{2}(r).
\]
If $j = 2$ then $r\le G = 2N^{\,2/3}$ for any $r\in E_{j}$, that is $N\ge r^{\,3/2}$ and $N^{2}+r^{2}\ll N^{2}$. Using the Corollary of Theorem 2.1 we get
\begin{multline*}
W_{r}\,\ll\,\pi(N)\biggl\{\biggl(\frac{r}{N}\biggr)^{\! 1/6}\mathcal{L}^{2}(\tau(r))^{2/3}\,+\,\mathcal{L}^{5/2}\,\frac{\tau(r)}{\sqrt[4\;]{r}}\biggr\}\,\ll\\
\ll\,\pi(N)\biggl\{N^{-1/18}\mathcal{L}^{2}(\tau(r))^{2/3}\,+\,\mathcal{L}^{5/2}\,\frac{\tau(r)}{\sqrt[4\;]{r}}\biggr\}.
\end{multline*}
Hence,
\begin{multline*}
I_{k}^{(2)}\,\ll\,\frac{1}{q}\sum\limits_{\substack{F<r\le G \\ r|q}}\pi^{k-2}(N)\biggl\{\biggl(\frac{\tau^{2}(r)\mathcal{L}^{6}}{N^{1/6}}\biggr)^{(k-2)/3}\,+\,\mathcal{L}^{5k/2-3}\,\frac{(\tau(r))^{k-2}}{r^{(k-2)/4}}\biggr\}N^{2}\tau^{2}(r)\,\ll\\
\ll\,\frac{\pi^{k}(N)}{q}\sum\limits_{\substack{F<r\le G \\ r|q}}\biggl\{\mathcal{L}^{2}\biggl(\frac{\tau^{2}(r)\mathcal{L}^{6}}{N^{1/6}}\biggr)^{(k-2)/3}\,+\,\mathcal{L}^{5k/2-1}\,\frac{(\tau(r))^{k-2}}{r^{(k-2)/4}}\biggr\}\tau^{2}(r)\,\ll
\end{multline*}
\begin{multline*}
\ll\,\frac{\pi^{k}(N)}{q}\biggl\{N^{-(k-2)/18}\mathcal{L}^{2(k-1)}\sum\limits_{r|q}(\tau(r))^{2(k+1)/3}\,+\,\mathcal{L}^{5k/2-1}\sum\limits_{\substack{F<r\le G \\ r|q}}\frac{(\tau(r))^{k}}{r^{(k-2)/4}}\biggr\}\,\ll\\
\ll\,\frac{\pi^{k}(N)}{q}\biggl\{N^{-(k-2)/18}\mathcal{L}^{2(k-1)}(\tau(q))^{(2k+5)/3}\,+\,\mathcal{L}^{5k/2-1}\sum\limits_{r>F}\frac{(\tau(r))^{k}}{r^{(k-2)/4}}\biggr\}.
\end{multline*}
The first term does not exceed $N^{-(k-2)/20}\le N^{-1/4}$. Next, in view of the inequality $k\ge 7$ we have $(k-2)/4\ge 5/4>1$ and the series over $r$ converges. By Mardzhanishvili's inequality \cite{Mardzhanishvili_1939}, we can easily conclude that this sum is bounded by
\[
\frac{(\ln{F})^{2^{k}-1}}{F^{(k-6)/4}}\,\ll\,\frac{(\ln{\mathcal{L}})^{2^{k}-1}}{\mathcal{L}^{C(k-6)/4}}.
\]
Thus we get
\begin{equation}\label{3-05}
I_{k}^{(2)}\,\ll\,\frac{\pi^{k}(N)}{q}\biggl\{N^{-1/4}\,+\,\mathcal{L}^{5k/2-1}\,\frac{(\ln{\mathcal{L}})^{2^{k}-1}}{\mathcal{L}^{C(k-6)/4}}\biggr\}\,\ll\,
\frac{\pi^{k}(N)}{q}\,\frac{(\ln{\mathcal{L}})^{2^{k}-1}}{\mathcal{L}^{A}},
\end{equation}
where
\[
A\,=\,\frac{C}{4}\,(k-6)\,-\,\frac{5k}{2}+1\,=\,\frac{29}{2}(k-7)\,+\,\frac{1}{2}.
\]
We apply the same arguments in the case $j = 3$, replacing the estimate of Corollary of Theorem 2.1 by the estimate of Lemma 3.1. Thus, for any fixed $\vep$ we have
\begin{multline*}
W_{r}\,\ll\,\bigl(N^{6/7}r^{3/28}\,+\,N^{32/35}r^{2/35}\bigr)r^{\vep/(4k)}\,\ll\,\pi(N)\bigl(N^{-1/7}q^{3/28}\,+\,N^{-3/35}q^{2/35}\bigr)q^{\vep/(3k)},\\
rJ_{r}(N)\,\ll\,(N^{2}+r^{2})\tau^{2}(r)\,\ll\,q^{2}\tau^{2}(q)\,\ll\,\pi^{2}(N)\,N^{-2}q^{2+\vep/6},
\end{multline*}
and hence
\begin{multline*}
I_{k}^{(3)}\,\ll\,\frac{\pi^{k}(N)}{q}\sum\limits_{\substack{G<r\le q \\ r|q}}\bigl(N^{-(k-2)/7}q^{3(k-2)/28}\,+\,N^{-3(k-2)/35}q^{2(k-2)/35}\bigr)N^{-2}q^{2}\,q^{\vep/2-2\vep/(3k)}\,\ll\\
\ll\,\frac{\pi^{k}(N)}{q}\,q^{0.5\vep}\bigl(N^{-(k-2)/7-2}q^{3(k-2)/28+2}\,+\,N^{-3(k-2)/35-2}q^{2(k-2)/35+2}\bigr)\,\ll\\
\ll\,\frac{\pi^{k}(N)}{q}\,q^{0.5\vep}\biggl\{\biggl(\frac{q^{\alpha_{k}}}{N}\biggr)^{(k+12)/7}\,+\,\biggl(\frac{q^{\beta_{k}}}{N}\biggr)^{(3k+64)/35}\biggr\},
\end{multline*}
where
\[
\alpha_{k}\,=\,\frac{3k+50}{4(k+12)},\quad \beta_{k}\,=\,\frac{2(k+33)}{3k+64}.
\]
Since $\max{(\alpha_{k},\beta_{k})} = c_{k}$, for $N\ge q^{c_{k}+\vep}$ we have
\begin{multline}\label{3-06}
I_{k}^{(3)}\,\ll\,\frac{\pi^{k}(N)}{q}\,q^{0.5\vep}\bigl(q^{-\vep(k+12)/7}\,+\,q^{-\vep(3k+64)/35}\bigr)\,\ll\\
\ll\,\frac{\pi^{k}(N)}{q}\,\bigl(q^{-\vep(2k+17)/14}\,+\,q^{-\vep(6k+93)/70}\bigr)\,\ll\,\frac{\pi^{k}(N)}{q}\,q^{-\vep}.
\end{multline}
Now the desired assertion follows from the relation (\ref{3-04}) and the estimates (\ref{3-05}), (\ref{3-06}).

The conditional proof follows the same lines, but we need to choose $F = N^{1/3}$, $G = N^{2/3}$. If the Generalized Riemann hypothesis is true then
\[
\pi(x;r,c)\,=\,\frac{\pi(x)}{\varphi(r)}\,+\,O(\sqrt{x}\ln{x})
\]
for any $(r,c)=1$ (see \cite[Ch. 20]{Davenport_1980}; the constant in $O$-symbol is absolute). Hence, passing to the calculation of $I_{k}^{(1)}$ and using the inequality (\ref{3-03}), we find
\begin{multline*}
W_{r}(af,bf;N)\,=\,\frac{\pi(N)}{\vf(r)}\,S(af,bf;r)\,+\,O\bigl(\vf(r)\mathcal{L}\sqrt{N}\bigr),\\
W_{\delta r}^{k}(a\delta f, b\delta f, N)\,=\,\frac{\pi^{k}(N)}{\vf^{k}(r)}\,S^{k}(af,bf;r)\,+\,O\biggl(\frac{\pi^{k-1}(N)}{\vf^{k-1}(r)}\,(\sqrt{r}\tau(r))^{k-1}\,\vf(r)\mathcal{L}\sqrt{N}\biggr)\,+\,\\
+\,O(\vf^{k}(r)N^{k/2}\mathcal{L}^{k})
\end{multline*}
and therefore
\begin{multline}\label{3-07}
I_{k}^{(1)}\,=\,\frac{\pi^{k}(N)}{q}\sum\limits_{\substack{1\le r\le F \\ r|q}}A_{k}(r)\,+\,O\biggl(\frac{N^{k/2}}{q}\,\mathcal{L}^{k}\sum\limits_{\substack{1\le r\le F \\ r|q}}\vf^{k+1}(r)\biggr)\,+\\
+\,O\biggl(\frac{\pi^{k-1}(N)}{q}\,\mathcal{L}\sqrt{N}\sum\limits_{\substack{1\le r\le F \\ r|q}}\frac{(\sqrt{r}\tau(r))^{k-1}}{\vf^{k-3}(r)}\biggr).
\end{multline}
The second term in (\ref{3-07}) does not exceed in order
\[
\frac{N^{k/2}\mathcal{L}^{k}}{q}\,F^{k+1}\tau(q)\,\ll\,\frac{\pi^{k}(N)}{q}\,\frac{F^{k+1}}{N^{k/2}}\,\mathcal{L}^{2k}\tau(q)\,\ll\,\frac{\pi^{k}(N)}{q}\,N^{-(k-2)/6}\mathcal{L}^{2k}\tau(q)\,\ll\,\frac{\pi^{k}(N)}{q}\,q^{-\vep}.
\]
Further, since
\[
\frac{r}{\vf(r)}\,\ll\,\ln\ln{r}\,\ll\,\ln{\mathcal{L}},
\]
then the last term in (\ref{3-07}) is estimated as
\begin{multline*}
\frac{\pi^{k}(N)}{q}\,\frac{\mathcal{L}^{2}}{\sqrt{N}}\,\sum\limits_{\substack{1\le r\le F \\ r|q}}\biggl(\frac{r}{\vf(r)}\biggr)^{k-3}\tau^{k-1}(q)\,\frac{r^{(k-1)/2}}{r^{k-3}}\,\ll\\
\ll\,\frac{\pi^{k}(N)}{q}\,\frac{\mathcal{L}^{2}}{\sqrt{N}}\,(\ln{\mathcal{L}})^{k-3}\tau^{k-1}(q)\sum\limits_{\substack{1\le r\le F \\ r|q}}r^{-(k-5)/2}\,\ll\\
\ll\,\frac{\pi^{k}(N)}{q}\,\frac{\mathcal{L}^{2}}{\sqrt{N}}\,(\ln{\mathcal{L}})^{k-3}\tau^{k}(q)\bigl(1\,+\,F^{-(k-5)/2}\bigr)\,\ll\\
\ll\,\frac{\pi^{k}(N)}{q}\,\mathcal{L}^{2}(\ln{\mathcal{L}})^{k-3}\tau^{k}(q)\bigl(N^{\,-1/2}\,+\,N^{-(k-2)/6}\bigr)\,\ll\,\frac{\pi^{k}(N)}{q}\,q^{-\vep}.
\end{multline*}
Next, the term $I_{k}^{(2)}$ is estimated by Corollary of Theorem 2.2:
\begin{multline}\label{3-08}
I_{k}^{(2)}\,\ll\,\frac{1}{q}\sum\limits_{\substack{F< r\le G \\ r|q}}\pi^{k-2}(N)\,\biggl\{\biggl(\frac{\sqrt{r}\tau(r)}{\vf(r)}\biggr)^{k-2}\,+\,\biggl(\frac{r}{N}\biggr)^{(k-2)/2}\mathcal{L}^{2(k-2)}\biggr\}\,N^{2}\tau^{2}(r)\,\ll\\
\ll\,\frac{\pi^{k}(N)}{q}\,\mathcal{L}^{2}\tau^{2}(q)\biggl\{\;\sum\limits_{\substack{F< r\le G \\ r|q}}\biggl(\frac{r}{\vf(r)}\biggr)^{k-2}\biggl(\frac{\tau(r)}{\sqrt{r}}\biggr)^{k-2}\,+\,\bigl(GN^{-1}\bigr)^{(k-2)/2}\mathcal{L}^{2(k-2)}\biggr\}\,\ll\\
\ll\,\frac{\pi^{k}(N)}{q}\,\mathcal{L}^{2}\tau^{2}(q)\biggl\{(\ln{\mathcal{L}})^{k-2}\sum\limits_{\substack{F< r\le G \\ r|q}}\biggl(\frac{\tau(q)}{\sqrt{r}}\biggr)^{k-2}\,+\,N^{-(k-2)/6}\mathcal{L}^{2(k-2)}\biggr\}\,\ll\\
\ll\,\frac{\pi^{k}(N)}{q}\,\mathcal{L}^{2}\bigl(F^{-(k-2)/2}(\ln{\mathcal{L}})^{k-2}\tau^{k+1}(q)\,+\,N^{-(k-2)/6}\tau^{2}(q)\mathcal{L}^{2(k-2)}\bigr)\,\ll\\
\ll\,\frac{\pi^{k}(N)}{q}\,N^{-(k-2)/6}\tau^{k+1}(q)\mathcal{L}^{2k-2}\,\ll\,\frac{\pi^{k}(N)}{q}\,q^{-\vep}.
\end{multline}
Finally, the estimate (\ref{3-06}) of $I_{k}^{(3)}$ in unconditional case is still true. Now the assertion follows from the relations (\ref{3-06})-(\ref{3-08}). $\Box$
\vspace{0.3cm}

\textsc{Remark.} Estimating the term $I_{k}^{(2)}$ unconditionally, we essentially use the fact that $k\ge 7$. In this case, the sum
\[
\sum\limits_{\substack{r|q \\ F<r\leqslant G}}\frac{(\tau(r))^{k}}{r^{(k-2)/4}}\,=\,\sum\limits_{\substack{r|q \\ r>F}}\frac{(\tau(r))^{k}}{r^{(k-2)/4}}\,+\,O\biggl(\frac{\tau^{k+1}(q)}{G}\biggr)\,=\,
\sum\limits_{\substack{r|q \\ r>F}}\frac{(\tau(r))^{k}}{r^{(k-2)/4}}\,+\,O\bigl(N^{-2/3+\varepsilon}\bigr)
\]
is small. It is naturally to ask whether this sum is small for $k\le 6$. It is easy to prove that this is true for almost all $q$. However, in general case, for $k = 6$, even the estimate
\[
\sum\limits_{\substack{r|q \\ r>F}}\frac{(\tau(r))^{k}}{r^{(k-2)/4}}\,=\,\sum\limits_{\substack{r|q \\ r>F}}\frac{\tau^{6}(r)}{r}\,=\,O(1)
\]
does not hold. To show this, it is sufficient to prove that for any fixed constant $A>1$ there exists an infinite sequence
 of moduli $q$ such that the sums
\[
\Sigma(q;A)\,=\,\sum\limits_{\substack{d|q \\ d>(\ln{q})^{A}}}\frac{1}{d}
\]
do not satisfy the estimate $\Sigma(q;A) = O(1)$, and, moreover $\Sigma(q;A)\to +\infty$ when $q\to +\infty$.

To do this, we need some definitions. Namely, let $\Psi(x,y)$ and $\Psi_{2}(x,y)$ denote the number of $n\le x$ (squarefree $n\le x$, respectively) that are free of prime factors $p\le y$ (here $2\le y\le x$). Then, given $\vep>0$, one has
\begin{align}\label{3-09}
& \Psi(x,y)\,=\,x\varrho(u)\biggl(1\,+\,O\Bigl(\frac{\ln{(1+u)}}{\ln{y}}\Bigr)\biggr),\\
\label{3-10}
& \Psi_{2}(x,y)\,=\,\frac{6}{\pi^{2}}\,\Psi(x,y)\biggl(1\,+\,O\Bigl(\frac{\ln^{2}{(u+1)}}{\ln{y}}\Bigr)\biggr),\quad u\,=\,\frac{\ln{x}}{\ln{y}}.
\end{align}
The first one holds uniformly in $e^{(\ln\ln{x})^{2+\vep}}\le y\le x$, and the second one holds uniformly in $x^{\vep}\le y\le x$
(for the proof of (\ref{3-09}), see, for example, \cite[Ch. III.5, Cor. 9.3]{Tenenbaum_1995}; the relation (\ref{3-10}) can be derived from the first one in the way pointed in \cite{Ivic_1995}). The symbol $\varrho(u)$ stands here for Dickman function defined by the following conditions:
\[
\varrho(u)\equiv 1\quad \text{when}\quad 0<u\le 1,\quad u\varrho'(u)+\varrho(u-1)\,=\,0\quad \text{for}\quad u>1.
\]
Now suppose that $A>1$ is arbitrary large but fixed number and suppose that $X\ge X_{0}(A)>2$. Denote $\displaystyle q = \prod\limits_{p\le X}p$ (here and below in this section, we use the same notation $q$ which stands for the modulus of the congruence; but this does not lead us to confusion). By prime number theorem,
\[
\ln{q}\,=\,\sum\limits_{p\le X}\ln{p}\,=\,X\,\bigl(1\,+\,O(e^{-c\sqrt{\ln{X}}})\bigr),\quad c>0,
\]
so we have $X^{A}= (\ln{q})^{A}(1+o(1))$. Next, by Mertens formula,
\begin{equation}\label{3-11}
S\,=\,\sum\limits_{d|q}\frac{1}{d}\,=\,\prod\limits_{p\le X}\biggl(1+\frac{1}{p}\biggr)\,=\,
\prod\limits_{p\le X}\biggl(1-\frac{1}{p^{2}}\biggr)\biggl(1-\frac{1}{p}\biggr)^{-1}\,=\,\frac{6}{\pi^{2}}\,e^{\gamma}(\ln{X})
\bigl(1\,+\,O(e^{-c\sqrt{\ln{X}}})\bigr).
\end{equation}
Further, let
\[
S_{1}\,=\,\sum\limits_{\substack{d|q \\ d\le X^{A}}}\frac{1}{d}.
\]
Then, using the notation $P^{+}(n)$ for the greatest prime divisor of $n\ge 2$, we get
\[
S_{1}\,=\,\sum\limits_{\substack{n\le X^{A} \\ P^{+}(n)\le X}}\frac{\mu^{2}(n)}{n}\,=\,\biggl(\,\sum\limits_{n\le X}\,+\,\sum\limits_{\substack{X<n\le X^{A} \\ P^{+}(n)\le X}}\,\biggr)\frac{\mu^{2}(n)}{n}\,=\,S_{2}+S_{3}.
\]
Basic property of M\"{o}bius function implies the formula
\[
S_{2}\,=\,\frac{6}{\pi^{2}}\,\Bigl(\ln{X}+\gamma+2\sum\limits_{p}\frac{\ln{p}}{p^{2}-1}\Bigr)\,+\,O\biggl(\frac{\ln{X}}{\sqrt{X}}\biggr).
\]
Next, by Abel summation formula we get
\begin{multline*}
S_{3}\,=\,\frac{1}{X^{A}}\bigl(\Psi_{2}(X^{A},X)-\Psi_{2}(X,X)\bigr)\,+\,\int_{X}^{X^{A}}\bigl(\Psi_{2}(u,X)-\Psi_{2}(X,X)\bigr)\frac{du}{u^{2}}\,=\\
=\,\frac{\Psi_{2}(X^{A},X)}{X^{A}}\,+\,\int_{X}^{X^{A}}\Psi_{2}(u,X)\,\frac{du}{u^{2}}\,-\,\frac{\Psi_{2}(X,X)}{X}.
\end{multline*}
Since
\[
\Psi_{2}(X,X)\,=\,\frac{6}{\pi^{2}}\,X\,+\,O\bigl(\sqrt{X}\bigr),
\]
then (\ref{3-10}) implies
\begin{multline*}
S_{3}\,=\,\frac{6}{\pi^{2}}\bigl(1+O(\Delta)\bigr)\int_{X}^{X^{A}}\varrho\biggl(\frac{\ln{u}}{\ln{X}}\biggr)\frac{du}{u}\,+\,
\frac{6}{\pi^{2}}\,\varrho(A)\bigl(1+O(\Delta)\bigr)\,-\,\frac{6}{\pi^{2}}\,+\,O\Bigl(\frac{1}{\sqrt{X}}\Bigr)\,=\\
=\,\frac{6}{\pi^{2}}\biggl(\,\int_{1}^{A}\varrho(v)dv\biggr)(\ln{X})\bigl(1+O(\Delta)\bigr),\quad\Delta = \frac{\ln^{2}(A+1)}{\ln{X}}.
\end{multline*}
Therefore,
\begin{equation}\label{3-12}
S_{1}\,=\,\frac{6}{\pi^{2}}\biggl(1\,+\,\int_{1}^{A}\varrho(v)dv\biggr)(\ln{X})\bigl(1+O(\Delta)\bigr)\,=\,
\frac{6}{\pi^{2}}\biggl(\,\int_{0}^{A}\varrho(v)dv\biggr)(\ln{X})\bigl(1+O(\Delta)\bigr).
\end{equation}
Finally, subtracting (\ref{3-12}) from (\ref{3-11}), we get
\[
\Sigma(q;A)\,=\,\sum\limits_{\substack{d|q \\ d>X^{A}}}\frac{1}{d}\,=\,\frac{6}{\pi^{2}}\biggl(e^{\gamma}\,-\,
\int_{0}^{A}\varrho(v)dv\biggr)(\ln{X})\bigl(1+O(\Delta)\bigr).
\]
Using the identity
\[
\int_{0}^{+\infty}\varrho(v)dv\,=\,2\varrho(2)+3\varrho(3)+4\varrho(4)+5\varrho(5)\ldots\,=\,e^{\gamma}
\]
(see \cite[Ch. III.5, Ex. 2]{Tenenbaum_1995}), we conclude that
\[
\Sigma(q;A)\,=\,c(A)(\ln{X})\bigl(1+O(\Delta)\bigr),\quad c(A)\,=\,\int_{A}^{+\infty}\varrho(v)dv\,>\,0.
\]
If $X\to +\infty$, then $\Sigma(q;A)\to +\infty$.
\vspace{0.3cm}

\textbf{4. Structure of the ``singular series''.}
\vspace{0.5cm}

In what follows, we study the properties of the ``singular series'' $\varkappa_{k}(q) = \varkappa_{k}(a,b,m;q)$. In this section, we give the ``probabilistic'' interpretation of the first term in the above formula for $I_{k}(N)$.

Since $\varkappa_{k}(q)$ is multiplicative over $q$, it is sufficient to study the case $q = p^{\,n}$ for prime $p$ and $n\ge 1$. Denote by $V_{k}(q) = V_{k}(a,b,m;q)$ the number of solutions of the congruence
\begin{equation}\label{4-01}
g(x_{1})\,+\,\ldots\,+\,g(x_{k})\,\equiv\,m\pmod{q}
\end{equation}
with the conditions $1\le x_{j}\le q$, $(x_{j},q)=1$, $j = 1,2,\ldots,k$. Obviously,
\[
V_{k}(q)\,=\,\frac{1}{q}\sum\limits_{f=1}^{q}e_{q}(-mf)S^{k}(af,bf;q).
\]
Then, for $n\ge 2$ we have
\begin{multline*}
A_{k}(p^{n})\,=\,\frac{1}{\varphi^{k}(p^{n})}\,\biggl\{\sum\limits_{f=1}^{p^{n}}e_{p^{n}}(-mf)S^{k}(af,bf;p^{n})\,-\,
\sum\limits_{f=1}^{p^{n-1}}e_{p^{n}}(-mfp)S^{k}(afp,bfp;p^{n})\biggr\}\,=\\
=\,\frac{1}{\varphi^{k}(p^{n})}\,\biggl\{\sum\limits_{f=1}^{p^{n}}e_{p^{n}}(-mf)S^{k}(af,bf;p^{n})\,-\,
p^{k}\sum\limits_{f=1}^{p^{n-1}}e_{p^{n-1}}(-mf)S^{k}(af,bf;p^{n-1})\biggr\}\,=\\
=\,\frac{1}{p^{k(n-1)}(p-1)^{k}}\,\bigl(p^{n}V_{k}(p^{n})-p^{k+n-1}V_{k}(p^{n-1})\bigr)\,=\,\frac{p}{(p-1)^{k}}\,\biggl(\,\frac{V_{k}(p^{n})}{p^{(k-1)(n-1)}}\,-\,\frac{V_{k}(p^{n-1})}{p^{(k-1)(n-2)}}\biggr).
\end{multline*}
By similar arguments, we find that
\[
A_{k}(p)\,=\,\frac{pV_{k}(p)-(p-1)^{k}}{(p-1)^{k}}\,=\,\frac{pV_{k}(p)}{(p-1)^{k}}\,-\,1.
\]
Thus we obtain
\begin{multline*}
\varkappa_{k}(p^{n})\,=\,\sum\limits_{\mu=0}^{n}A_{k}(p^{\mu})\,=\,1\,+\,\biggl(\frac{pV_{k}(p)}{(p-1)^{k}}\,-\,1\biggr)\,+\,
\sum\limits_{\mu=2}^{n}\frac{p}{(p-1)^{k}}\,\biggl(\,\frac{V_{k}(p^{\,\mu})}{p^{(k-1)(\mu-1)}}\,-\,\frac{V_{k}(p^{\,\mu-1})}{p^{(k-1)(\mu-2)}}\biggr)\,=\\
=\,\frac{p}{(p-1)^{k}}\cdot \frac{V_{k}(p^{\,n})}{p^{(k-1)(n-1)}}\,=\,\frac{p^{\,n}V_{k}(p^{\,n})}{\vf^{k}(p^{\,n})}.
\end{multline*}
Hence, in the case of an arbitrary modulus $q$ we have
\[
\varkappa_{k}(q)\,=\,\frac{q}{\varphi^{k}(q)}\prod\limits_{p^{\,n}||q}V_{k}(p^{\,n})\,=\,\frac{qV_{k}(q)}{\varphi^{k}(q)}
\]
(here we use the multiplicativity of $V_{k}(q)$ over $q$; this fact easily follows from Chinese remainder theorem).
Now the first term in the expression for $I_{k}(N)$ can be expressed as follows:
\begin{equation}\label{4-02}
\frac{\pi^{k}(N)}{q}\,\frac{qV_{k}(q)}{\varphi^{k}(q)}\,=\,\biggl(\frac{\pi(N)}{\varphi(q)}\biggr)^{\!k}\,V_{k}(q)
\end{equation}
The ratio $\pi(N)/\varphi(q)$ is the density of primes $p\le N$ in the reduced residual system ($N<q$). Therefore, the coefficient
$(\pi(N)/\varphi(q))^{k}$ in (\ref{4-02}) is the ``probability''  of that all the components of the tuple $(x_{1},\ldots,x_{k})$ satisfying (\ref{4-01}), are prime numbers. In other words, this is the probability of that the solution of (\ref{4-01}) is also the solution of the initial congruence (\ref{1-04}).
\vspace{0.5cm}

\textbf{5. Explicit formulae for the quantities $\boldsymbol{A_{k}(p^{n})}$, $\boldsymbol{3\le k\le 5}$, via Ramanujan sums.}
\vspace{0.5cm}

To study the properties of $\varkappa_{k}(q)$, we need the explicit expressions for the quantities $S(a,b;p^{n})$, $A_{k}(a,b,m;p^{n})$, $n\ge 2$.
\vspace{0.3cm}

\textsc{Lemma 5.1.} \emph{Let $p\ge 3$ and $(ab,p)=1$. If $ab$ is a quadratic non-residue modulo $p$ then $S(a,b;p^{n}) = 0$ for any $n\ge 2$. Otherwise, setting $q$ for $p^{n}$ and $\nu$ for any solution of the congruence $ab\equiv \nu^{2}\pmod{q}$, we have}
\begin{equation*}
S(a,b;q)\,=\,S(\nu,\nu;q)\,=
\begin{cases}
2\sqrt{q}\cos{\Bigl(\frac{\displaystyle 4\pi\nu}{\displaystyle q}\Bigr)} & \textit{for even } n,\\
2\sqrt{q}\Bigl(\frac{\displaystyle \nu}{\displaystyle q}\Bigr)\cos{\Bigl(\frac{\displaystyle 4\mathstrut\pi\nu}{\displaystyle q}+\frac{\displaystyle \pi s}{\displaystyle 2\mathstrut}\Bigr)} & \textit{for odd } n,
\end{cases}
\end{equation*}
\emph{where $s = 0$ for $p\equiv 1\pmod{4}$ and $s = 1$ for} $p\equiv 3\pmod{4}$.
\vspace{0.3cm}

For the proof, see \cite{Salie_1931}.
\vspace{0.3cm}

\textsc{Corollary.} \emph{Under the conditions of Lemma} 5.1, \emph{for any $n\ge 2$ and $q = p^{n}$ we have}
\[
|S(a,b;q)|\,<\,2\sqrt{q}.
\]

\textsc{Remark.} For the most part of number-theoretic applications, only the Corollary of Lemma 5.1 is sufficient. It is surprisingly enough that the below explicit formulas for $A_{k}(q)$ are used essentially in the deriving of
the precise expressions for complete Kloosterman sums given by Lemma 5.1. The connection between Lemma 5.1 and the expressions for $A_{k}(q)$ is given by well-known Ramanujan sums.
\vspace{0.3cm}

Recall that Ramanujan sum $c_{q}(a)$ (see \cite{Ramanujan_1918}) is defined by the following relations:
\[
c_{q}(a)\,=\,\mathop{{\sum}'}\limits_{f=1}^{q}e^{2\pi i\frac{\scriptstyle af}{\scriptstyle q}}\,=\,\mathop{{\sum}'}\limits_{f=1}^{q}\cos{\Bigl(\frac{2\pi af}{q}\Bigr)}.
\]
It is well-known that
\[
c_{q}(a)\,=\,\frac{\varphi(q)}{\varphi(q/\delta)}\,\mu\Bigl(\frac{q}{\delta}\Bigr),\quad \delta = (q,a),
\]
(see, for example, \cite[Ch. XVI, Th. 272]{Hardy_Wright_1975}). In particular, if $q = p^{n}$, $n\ge 2$, and $\delta= p^{r}$ then
\begin{equation}\label{7-0}
c_{q}(a)\,=\,
\begin{cases}
0, & \text{if}\;\;0\le r\le n-2,\\
-p^{n-1}, & \text{if}\;\; r= n-1,\\
\vf(p^{n}), & \text{if}\;\;  r = n.
\end{cases}
\end{equation}
We also introduce the quantities
\[
u_{q}(a)\,=\,\RRe{w_{q}(a)},\quad v_{q}(a)\,=\,\IIm{w_{q}(a)},\quad\text{where}\quad w_{q}(a)\,=\,\sum\limits_{f=1}^{q}\Bigl(\frac{f}{q}\Bigr)e^{2\pi i\frac{\scriptstyle af}{\scriptstyle q}}
\]
and $\displaystyle \Bigl(\frac{f}{q}\Bigr)$ stands for Jacobi symbol in this section.
\vspace{0.3cm}

\textsc{Lemma 5.2.} \emph{If $p\ge 3$ is prime, $n$  is odd and $q = p^{n}$ then the following relations hold:}
\begin{multline*}
u_{q}(a)\,=\,
\begin{cases}
\displaystyle \Bigl(\frac{b}{p}\Bigr)p^{\,n-1/2}, & \textit{when}\;\;a=bp^{n-1},\; p\equiv 1\pmod{4},\\
0, & \textit{otherwise},
\end{cases}\\
v_{q}(a)\,=\,
\begin{cases}
\displaystyle \Bigl(\frac{b}{p}\Bigr)p^{\,n-1/2}, & \textit{when}\;\;a=bp^{n-1},\; p\equiv 3\pmod{4},\\
0, & \textit{otherwise}.
\end{cases}
\end{multline*}
\textsc{Proof.} In fact, let $a = bp^{r}$ where $(b,p)=1$ and $0\le r\le n-1$. Then, setting $f = g+p^{n-r}h$, where $1\le g\le p^{n-r}$ and $1\le h\le p^{r}$, we get
\begin{multline*}
w_{q}(a)\,=\,\sum\limits_{g=1}^{p^{n-r}}\sum\limits_{h=1}^{p^{r}}\Bigl(\frac{g}{p}\Bigr)\exp{\biggl(2\pi i\frac{b}{p^{\,n-r}}\,(g+p^{n-r}h)\biggr)}\,=\,p^{r}
\sum\limits_{g=1}^{p^{n-r}}\Bigl(\frac{g}{p}\Bigr)\exp{\biggl(2\pi i\frac{bg}{p^{\,n-r}}\biggr)}\,=\\
=\,p^{r}\Bigl(\frac{b}{p}\Bigr)\sum\limits_{g=1}^{p^{n-r}}\Bigl(\frac{g}{p}\Bigr)\exp{\biggl(2\pi i\frac{g}{p^{\,n-r}}\biggr)}.
\end{multline*}
If $n-r=1$ then the sum over $g$ coincides with Gaussian sum
\begin{equation*}
\sum\limits_{g=1}^{p}\Bigl(\frac{g}{p}\Bigr)e^{2\pi i\frac{\scriptstyle g}{\scriptstyle p}}\,=\,
\begin{cases}
\sqrt{p}, & p\equiv 1\pmod{4},\\
i\sqrt{p}, & p\equiv 3\pmod{4}.
\end{cases}
\end{equation*}
Otherwise, setting $g = s+p^{n-r-1}t$ where $1\le s\le p^{n-r-1}$, $1\le t\le p$, we find
\[
w_{q}(a)\,=\,p^{\,r}\Bigl(\frac{b}{p}\Bigr)\sum\limits_{s=1}^{p^{n-r-1}}\Bigl(\frac{s}{p}\Bigr)\exp{\biggl(2\pi i\frac{s}{p^{\,n-r}}\biggr)}\sum\limits_{t=1}^{p}
e^{2\pi i\frac{\scriptstyle t}{\scriptstyle p}}\,=\,0.
\]
Lemma is proved. $\Box$
\vspace{0.3cm}

\textsc{Lemma 5.3.} \emph{Suppose that $p\ge 3$, $n\ge 2$, $q = p^{n}$, and let $a,b,m$ be any triple satisfying to the conditions: $1\le a,b,m\le q$, $(ab,p)=1$,
$\displaystyle \Bigl(\frac{ab}{p}\Bigr) = 1$. Finally, let $\nu$ be any solution of the congruence $ab\equiv\nu^{2}\pmod{q}$. Then the following relations hold:}
\begin{multline}\label{7-1}
A_{3}(q)\,=\,
\begin{cases}
& \displaystyle\frac{q\sqrt{q}}{\vf^{3}(q)}\,\bigl(c_{q}(m+6\nu)+c_{q}(m-6\nu)+3c_{q}(m+2\nu)+3c_{q}(m-2\nu)\bigr), \\
& \qquad\qquad \textit{when}\;\;n\;\;\textit{is even},\\
& \displaystyle\frac{q\sqrt{q}}{\vf^{3}(q)}\Bigl(\frac{\nu}{p}\Bigr)\,\bigl(u_{q}(m+6\nu)+u_{q}(m-6\nu)+3u_{q}(m+2\nu)+3u_{q}(m-2\nu)\bigr), \\
& \qquad\qquad \textit{when}\;\;n\;\;\textit{is odd and}\;\;p\equiv 1\pmod{4},\\
& \displaystyle\frac{q\sqrt{q}}{\vf^{3}(q)}\Bigl(\frac{\nu}{p}\Bigr)\,\bigl(v_{q}(m+6\nu)-v_{q}(m-6\nu)-3v_{q}(m+2\nu)+3v_{q}(m-2\nu)\bigr), \\
& \qquad\qquad \textit{when}\;\;n\;\;\textit{is odd and}\;\;p\equiv 3\pmod{4}.
\end{cases}
\end{multline}
\textsc{Proof.} By Lemma 5.1, for even $n$ we have
\[
S(af,bf;q)\,=\,S(\nu f, \nu f;q)\,=\,2\sqrt{q}\cos{\Bigl(\frac{4\pi\nu f}{q}\Bigr)}.
\]
Hence,
\begin{multline*}
A_{3}(q)\,=\,\frac{1}{\varphi^{3}(q)}\mathop{{\sum}'}\limits_{f=1}^{q}\cos{\Bigl(\frac{2\pi mf}{q}\Bigr)}\biggl(2\sqrt{q}\cos{\Bigl(\frac{4\pi\nu f}{q}\Bigr)}\biggr)^{3}
\,=\\
=\,\frac{8q\sqrt{q}}{\varphi^{3}(q)}\mathop{{\sum}'}\limits_{f=1}^{q}\cos{\Bigl(\frac{2\pi mf}{q}\Bigr)}\cos^{3}{\Bigl(\frac{4\pi\nu f}{q}\Bigr)}.
\end{multline*}
Next, using the identity $\displaystyle\cos^{3}{\vartheta} = \frac{1}{4}\cos{3\vartheta}+\frac{3}{4}\cos{\vartheta}$, we obtain
\begin{multline*}
A_{3}(q)\,=\,\frac{8q\sqrt{q}}{\varphi^{3}(q)}\mathop{{\sum}'}\limits_{f=1}^{q}\cos{\Bigl(\frac{2\pi mf}{q}\Bigr)}\biggl(\frac{1}{4}\cos{\Bigl(\frac{12\pi\nu f}{q}\Bigr)}\,+\,\frac{3}{4}\cos{\Bigl(\frac{4\pi\nu f}{q}\Bigr)}\biggr)\,=\\
=\,\frac{2q\sqrt{q}}{\varphi^{3}(q)}\mathop{{\sum}'}\limits_{f=1}^{q}\biggl(\cos{\Bigl(\frac{2\pi m f}{q}\Bigr)}\cos{\Bigl(\frac{12\pi\nu f}{q}\Bigr)}\,+\,
3\cos{\Bigl(\frac{2\pi m f}{q}\Bigr)}\cos{\Bigl(\frac{4\pi\nu f}{q}\Bigr)}\biggr)\,=\\
=\,\frac{q\sqrt{q}}{\varphi^{3}(q)}\mathop{{\sum}'}\limits_{f=1}^{q}\biggl\{\cos{\Bigl(\frac{2\pi f}{q}(m+6\nu)\Bigr)}
\,+\,\cos{\Bigl(\frac{2\pi f}{q}(m-6\nu)\Bigr)}\,+\\
+\,3\cos{\Bigl(\frac{2\pi f}{q}(m+2\nu)\Bigr)}+3\cos{\Bigl(\frac{2\pi f}{q}(m-2\nu)\Bigr)}\biggr\}\,=\\
=\,\frac{q\sqrt{q}}{\varphi^{3}(q)}\bigl(c_{q}(m+6\nu)+c_{q}(m-6\nu)+3c_{q}(m+2\nu)+3c_{q}(m-2\nu)\bigr).
\end{multline*}
In the case of odd $n$ and $p\equiv 1\pmod{4}$, Lemma 5.1 implies:
\begin{multline*}
A_{3}(q)\,=\,\frac{1}{\varphi^{3}(q)}\mathop{{\sum}'}\limits_{f=1}^{q}\cos{\Bigl(\frac{2\pi mf}{q}\Bigr)}\biggl(2\Bigl(\frac{\nu f}{q}\Bigr)\sqrt{q}\cos{\Bigl(\frac{4\pi\nu f}{q}\Bigr)}\biggr)^{3}\,=\\
=\,\frac{8q\sqrt{q}}{\varphi^{3}(q)}\Bigl(\frac{\nu}{q}\Bigr)\mathop{{\sum}'}\limits_{f=1}^{q}\Bigl(\frac{f}{q}\Bigr)
\cos{\Bigl(\frac{2\pi m f}{q}\Bigr)}\cos^{3}{\Bigl(\frac{4\pi \nu f}{q}\Bigr)}\,=\\
=\,\frac{8q\sqrt{q}}{\varphi^{3}(q)}\Bigl(\frac{\nu}{q}\Bigr)\mathop{{\sum}'}\limits_{f=1}^{q}\Bigl(\frac{f}{q}\Bigr)
\cos{\Bigl(\frac{2\pi m f}{q}\Bigr)}\biggl\{\frac{1}{4}\cos{\Bigl(\frac{12\pi \nu f}{q}\Bigr)}\,+\,
\frac{3}{4}\cos{\Bigl(\frac{4\pi \nu f}{q}\Bigr)}\biggr\}\,=\\
=\,\frac{2q\sqrt{q}}{\varphi^{3}(q)}\Bigl(\frac{\nu}{q}\Bigr)\mathop{{\sum}'}\limits_{f=1}^{q}\Bigl(\frac{f}{q}\Bigr)
\biggl\{\cos{\Bigl(\frac{2\pi m f}{q}\Bigr)}\cos{\Bigl(\frac{12\pi \nu f}{q}\Bigr)}\,+\,3\cos{\Bigl(\frac{2\pi m f}{q}\Bigr)}\cos{\Bigl(\frac{4\pi \nu f}{q}\Bigr)}\biggr\}\,=
\end{multline*}
\begin{multline*}
=\,\frac{q\sqrt{q}}{\varphi^{3}(q)}\Bigl(\frac{\nu}{q}\Bigr)\mathop{{\sum}'}\limits_{f=1}^{q}\Bigl(\frac{f}{q}\Bigr)
\biggl\{\cos{\Bigl(\frac{2\pi}{q}(m+6\nu)\Bigr)}\,+\,\cos{\Bigl(\frac{2\pi}{q}(m-6\nu)\Bigr)}\,+\\
+\,3\cos{\Bigl(\frac{2\pi}{q}(m+2\nu)\Bigr)}\,+\,3\cos{\Bigl(\frac{2\pi}{q}(m-2\nu)\Bigr)}\biggr\}\,=\\
=\,\frac{q\sqrt{q}}{\varphi^{3}(q)}\Bigl(\frac{\nu}{q}\Bigr)\bigl(u_{q}(m+6\nu)+u_{q}(m-6\nu)+3u_{q}(m+2\nu)+3u_{q}(m-2\nu)\bigr).
\end{multline*}
In the case of odd $n$ and $p\equiv 3\pmod{4}$, the identity $\displaystyle \sin^{3}{\vartheta} = -\frac{1}{4}\sin{3\vartheta}+\frac{3}{4}\sin{\vartheta}$ together with Lemma 5.1 imply:
\begin{multline*}
A_{3}(q)\,=\,-\,\frac{8q\sqrt{q}}{\varphi^{3}(q)}\Bigl(\frac{\nu}{q}\Bigr)\mathop{{\sum}'}\limits_{f=1}^{q}\Bigl(\frac{f}{q}\Bigr)
\cos{\Bigl(\frac{2\pi mf}{q}\Bigr)}\sin^{3}{\Bigl(\frac{4\pi\nu f}{q}\Bigr)}\,=\\
=\,\frac{8q\sqrt{q}}{\varphi^{3}(q)}\Bigl(\frac{\nu}{q}\Bigr)\mathop{{\sum}'}\limits_{f=1}^{q}\Bigl(\frac{f}{q}\Bigr)
\cos{\Bigl(\frac{2\pi mf}{q}\Bigr)}\biggl\{\frac{1}{4}\sin{\Bigl(\frac{12\pi\nu f}{q}\Bigr)}\,-\,\frac{3}{4}
\sin{\Bigl(\frac{4\pi\nu f}{q}\Bigr)}\biggr\}\,=\\
=\,\frac{2q\sqrt{q}}{\varphi^{3}(q)}\Bigl(\frac{\nu}{q}\Bigr)\mathop{{\sum}'}\limits_{f=1}^{q}\Bigl(\frac{f}{q}\Bigr)\biggl\{
\cos{\Bigl(\frac{2\pi m f}{q}\Bigr)}\sin{\Bigl(\frac{12\pi\nu f}{q}\Bigr)}\,-\,3\cos{\Bigl(\frac{2\pi m f}{q}\Bigr)}\sin{\Bigl(\frac{4\pi\nu f}{q}\Bigr)}\biggr\}\,=\\
=\,\frac{q\sqrt{q}}{\varphi^{3}(q)}\Bigl(\frac{\nu}{q}\Bigr)\mathop{{\sum}'}\limits_{f=1}^{q}\Bigl(\frac{f}{q}\Bigr)
\biggl\{\sin{\Bigl(\frac{2\pi f}{q}(m+6\nu)\Bigr)}\,-\,\sin{\Bigl(\frac{2\pi f}{q}(m-6\nu)\Bigr)}\,-\\
-\,3\sin{\Bigl(\frac{2\pi f}{q}(m+2\nu)\Bigr)}\,+\,3\sin{\Bigl(\frac{2\pi f}{q}(m-2\nu)\Bigr)}\biggr\}\,=\\
=\,\frac{q\sqrt{q}}{\varphi^{3}(q)}\Bigl(\frac{\nu}{q}\Bigr)\bigl(v_{q}(m+6\nu)-v_{q}(m-6\nu)-3v_{q}(m+2\nu)+3v_{q}(m-2\nu)\bigr).
\end{multline*}
Lemma is proved. $\Box$
\vspace{0.3cm}

\textsc{Lemma 5.4.} \emph{Under the conditions of Lemma} 5.2, \emph{we have}
\begin{equation}\label{7-2}
A_{4}(q)\,=\,\frac{q^{2}}{\varphi^{4}(q)}\Bigl(c_{q}(m+8\nu)+c_{q}(m-8\nu)\,+\,4(-1)^{\frac{\scriptstyle n}{\scriptstyle 2\mathstrut}(p-1)}\bigl(c_{q}(m+4\nu)+c_{q}(m-4\nu)\bigr)\,+\,6c_{q}(m)\Bigr)
\end{equation}

\textsc{Remark.} Obviously, $(-1)^{\frac{\scriptstyle n}{\scriptstyle 2\mathstrut}(p-1)}=-1$ only in the case when $p\equiv 3\pmod{4}$, $n\equiv 1\pmod{2}$.
\vspace{0.3cm}

\textsc{Proof.} For even $n$, Lemma 5.1 implies
\[
A_{4}(q)\,=\,\frac{16q^{2}}{\varphi^{4}(q)}\mathop{{\sum}'}\limits_{f=1}^{q}\cos{\Bigl(\frac{2\pi m f}{q}\Bigr)}
\cos^{4}{\Bigl(\frac{4\pi \nu f}{q}\Bigr)}.
\]
Hence, using the identity $\displaystyle \cos^{4}{\vartheta}=\frac{1}{8}\cos{4\vartheta}+\frac{1}{2}\cos{2\vartheta}+\frac{3}{8}$, we get
\begin{multline*}
A_{4}(q)\,=\,\frac{2q^{2}}{\varphi^{4}(q)}\mathop{{\sum}'}\limits_{f=1}^{q}\cos{\Bigl(\frac{2\pi m f}{q}\Bigr)}\biggl\{
\cos{\Bigl(\frac{16\pi \nu f}{q}\Bigr)}+4\cos{\Bigl(\frac{8\pi \nu f}{q}\Bigr)}+3\biggr\}\,=\\
=\,\frac{q^{2}}{\varphi^{4}(q)}\mathop{{\sum}'}\limits_{f=1}^{q}\biggl\{\cos{\Bigl(\frac{2\pi f}{q}(m+8\nu)\Bigr)}\,+\,
\cos{\Bigl(\frac{2\pi f}{q}(m-8\nu)\Bigr)}\,+\\
+\,4\cos{\Bigl(\frac{2\pi f}{q}(m+4\nu)\Bigr)}\,+\,4\cos{\Bigl(\frac{2\pi f}{q}(m-4\nu)\Bigr)}\,+\,6\cos{\Bigl(\frac{2\pi fm}{q}\Bigr)}\biggr\}\,=\\
=\,\frac{q^{2}}{\varphi^{4}(q)}\bigl(c_{q}(m+8\nu)+c_{q}(m-8\nu)\,+\,4c_{q}(m+4\nu)+4c_{q}(m-4\nu)+6c_{q}(m)\bigr).
\end{multline*}
If $n$ is odd and $p\equiv 1\pmod{4}$ then Lemma 5.1 yields
\[
S^{4}(af,bf;q)\,=\,S^{4}(\nu f,\nu f; q)\,=\,\biggl\{2\Bigl(\frac{\nu f}{q}\Bigr)\sqrt{q}\cos{\Bigl(\frac{4\pi\nu f}{q}\Bigr)}\biggr\}^{4}\,=\,
16q^{2}\cos^{4}{\Bigl(\frac{4\pi\nu f}{q}\Bigr)}.
\]
Therefore, $A_{4}(q)$ coincides with the expression given above. Finally, if $n$ is odd and $p\equiv 3\pmod{4}$ then, by Lemma 5.1,
\[
S^{4}(af,bf;q)\,=\,\biggl\{-2\Bigl(\frac{\nu f}{q}\Bigr)\sqrt{q}\sin{\Bigl(\frac{4\pi\nu f}{q}\Bigr)}\biggr\}^{4}\,=\,
16q^{2}\sin^{4}{\Bigl(\frac{4\pi\nu f}{q}\Bigr)}.
\]
Thus, using the identity $\displaystyle \sin^{4}{\vartheta}=\frac{1}{8}\cos{4\vartheta}-\frac{1}{2}\cos{2\vartheta}+\frac{3}{8}$, we find
\begin{multline*}
A_{4}(q)\,=\,\frac{2q^{2}}{\varphi^{4}(q)}\mathop{{\sum}'}\limits_{f=1}^{q}\cos{\Bigl(\frac{2\pi m f}{q}\Bigr)}\biggl\{
\cos{\Bigl(\frac{16\pi \nu f}{q}\Bigr)}-4\cos{\Bigl(\frac{8\pi \nu f}{q}\Bigr)}+3\biggr\}\,=\\
=\,\frac{q^{2}}{\varphi^{4}(q)}\bigl(c_{q}(m+8\nu)+c_{q}(m-8\nu)\,-\,4c_{q}(m+4\nu)-4c_{q}(m-4\nu)+6c_{q}(m)\bigr).
\end{multline*}
Lemma is proved. $\Box$
\vspace{0.3cm}

\textsc{Lemma 5.5.} \emph{Under the conditions of Lemma} 5.2, \emph{the following relations hold:}
\begin{multline*}
A_{5}(q)\,=\,
\begin{cases}
& \displaystyle\frac{q^{2}\sqrt{q}}{\vf^{5}(q)}\,\bigl(c_{q}(m+10\nu)+c_{q}(m-10\nu)+5c_{q}(m+6\nu)+5c_{q}(m-6\nu)+\\
& \quad\quad +\,10c_{q}(m+2\nu)+10c_{q}(m-2\nu)\bigr), \quad \textit{when}\;\;n\;\;\textit{is even},\\
& \displaystyle\frac{q^{2}\sqrt{q}}{\vf^{5}(q)}\Bigl(\frac{\nu}{p}\Bigr)\,\bigl(u_{q}(m+10\nu)+u_{q}(m-10\nu)+5u_{q}(m+6\nu)+5u_{q}(m-6\nu)+\\
& \quad\quad +\,10u_{q}(m+2\nu)+10u_{q}(m-2\nu)\bigr), \quad \textit{when}\;\;n\;\;\textit{is odd and}\;\;p\equiv 1\pmod{4},\\
& \displaystyle\frac{q^{2}\sqrt{q}}{\vf^{5}(q)}\Bigl(\frac{\nu}{p}\Bigr)\,\bigl(-v_{q}(m+10\nu)+v_{q}(m-10\nu)+5v_{q}(m+6\nu)-5v_{q}(m-6\nu)-\\
& \quad\quad -\,10v_{q}(m+2\nu)+10v_{q}(m-2\nu)\bigr), \quad \textit{when}\;\;n\;\;\textit{is odd and}\;\;p\equiv 3\pmod{4}.
\end{cases}
\end{multline*}

\textsc{Proof.} The proof follows the same lines as above and uses the formulas
\[
\cos^{5}\vartheta\,=\,\frac{1}{2^{4\mathstrut}}\bigl(\cos{5\vartheta}+5\cos{3\vartheta}+10\cos{\vartheta}\bigr),\quad
\sin^{5}\vartheta\,=\,\frac{1}{2^{4\mathstrut}}\bigl(\sin{5\vartheta}-5\sin{3\vartheta}+10\sin{\vartheta}\bigr).\quad \Box
\]
\vspace{0.3cm}

\textsc{Remark.} Using the same arguments, on can obtain the expressions for $A_{k}(q)$ in the case of arbitrary $k$. However, Lemmas 5.3-5.5 are sufficient for our purposes.
\vspace{0.5cm}

\textbf{6. Some estimates for the quantities $\boldsymbol{A_{k}(p^{n})}$.}
\vspace{0.5cm}

Explicit formulas for the sums $A_{k}(p^{n})$ given in \S\,5 together with the Corollary of Lemma 5.1 are used here to estimate absolute values of these quantities.
\vspace{0.3cm}

\textsc{Lemma 6.1.} \emph{Suppose that $p\ge 3$ is prime. Then for any $s\ge n\ge 2$, $k\ge 5$ the following inequality holds:}
\[
\Bigl|\sum\limits_{r = n}^{s}A_{k}(p^{\,r})\Bigr|\,<\,2\biggl(\frac{2p}{p-1}\biggr)^{k-1}\frac{p^{-n(k/2-1)}}{1-p^{1-k/2}}.
\]
\textsc{Proof.} Lemma 5.1 implies the bound
\[
|A_{k}(p^{r})|\,<\,2\Bigl(\frac{2p}{p-1}\Bigr)^{\!k-1}p^{-\,r(k/2-1)}.
\]
Summing this estimate over $n\le r\le s$, we get the desired result. $\Box$
\vspace{0.3cm}

\textsc{Lemma 6.2.} \emph{Suppose that $p\ge 3$ is prime. Then the inequality}
\[
|A_{k}(p)|\,<\,\frac{(2\sqrt{p})^{k}}{(p-1)^{k-1}}\biggl(\frac{1}{4}\,+\,\frac{1}{4p}\biggr)
\]
\emph{holds for any $k\ge 5$.}
\vspace{0.3cm}

\textsc{Proof.} By Lemma 2.1, $|S(u,v;p)|<2\sqrt{p}$ for any $u,v$ such that $(uv,p)=1$. Therefore,
\begin{multline*}
|A_{k}(p)|\,\le\,\frac{1}{\vf^{k}(p)}\sum\limits_{f=1}^{p-1}|S(af,bf;p)|^{k}\,<\,\frac{(2\sqrt{p})^{k-2}}{\vf^{k}(p)}\sum\limits_{f=1}^{p-1}|S(af,bf;p)|^{2}\,=\\
=\,\frac{(2\sqrt{p})^{k-2}}{\vf^{k}(p)}\biggl(\,\sum\limits_{f=1}^{p}|S(af,bf;p)|^{2}\,-\,|S(0,0;p)|^{2}\biggr)\,=\,\frac{(2\sqrt{p})^{k-2}}{\vf^{k}(p)}\bigl(pI-(p-1)^{2}\bigr),
\end{multline*}
where $I$ denotes the number of solutions of the congruence
\[
g(x)\,\equiv\,g(y)\pmod{p},\quad 1\le x,y\le p-1.
\]
Obviously, $I\le 2(p-1)$. Hence,
\[
|A_{k}(p)|\,\le\,\frac{(2\sqrt{p})^{k-2}}{(p-1)^{k}}\bigl(2p(p-1)-(p-1)^{2}\bigr)\,=\,\frac{(2\sqrt{p})^{k-2}}{(p-1)^{k-1}}\,(p+1)\,=\,\frac{(2\sqrt{p})^{k}}{(p-1)^{k-1}}\biggl(\frac{1}{4}+\frac{1}{4p}\biggr).
\]
Lemma is proved. $\Box$
\vspace{0.3cm}

\textsc{Lemma 6.3.} \emph{Let $p\ge 3$ be a prime and suppose that $n\ge 3$ when $p=3$ and $n\ge 2$ when $p\ge 5$. Then, for any $s\ge n$ and for any $a,b$ satisfying the condition $\displaystyle \Bigl(\frac{ab}{p}\Bigr)=1$ the following inequality holds:}
\begin{align*}
&\biggl|\,\sum\limits_{r=n}^{s}A_{3}(p^{r})\biggr|\,\le\,
\begin{cases}
\displaystyle 3\Bigl(\frac{p}{p-1}\Bigr)^{4}p^{-\,\frac{\scriptstyle n}{\scriptstyle 2\mathstrut}}, &\textit{when}\quad n\equiv 0\pmod{2},\\
\displaystyle 3\Bigl(2-\frac{1}{p}\Bigr)\Bigl(\frac{p}{p-1}\Bigr)^{4}p^{-\,\frac{\scriptstyle n+1}{\scriptstyle 2\mathstrut}}, &\textit{when}\quad n\equiv 1\pmod{2};\\
\end{cases}
\end{align*}
\emph{In particular,}
\begin{align*}
& \biggl|\,\sum\limits_{r=2}^{s}A_{3}(p^{r})\biggr|\,\le\,\frac{3p^{3}}{(p-1)^{4}}, && \biggl|\,\sum\limits_{r=4}^{s}A_{3}(p^{r})\biggr|\,\le\,\frac{3p^{2}}{(p-1)^{4}}, \\
& \biggl|\,\sum\limits_{r=3}^{s}A_{3}(p^{r})\biggr|\,\le\,\frac{3p(2p-1)}{(p-1)^{4}}, && \biggl|\,\sum\limits_{r=5}^{s}A_{3}(p^{r})\biggr|\,\le\,\frac{3(2p-1)}{(p-1)^{4}}.
\end{align*}

\textsc{Proof.} First we prove that each of the formulas (\ref{7-1}) for $A_{3}(q)$ contains at most one non-zero term. Indeed, in the opposite case, the relations (\ref{7-0}) together with Lemma 5.2 imply that there exist two numbers among $m\pm 6\nu$, $m\pm 2\nu$ dividing by $p^{n-1}$. Hence, the difference of such numbers is also divisible by $p^{n-1}$. But this difference has the form
$\pm 4\nu$, $\pm 8\nu$, $\pm 12\nu$ and can not be divided by $p^{n-1}$. This leads us to contradiction.

Next, let $n\le r\le s$. Since $|c_{q}(a)|\le \varphi(q)$, $|u_{q}(a)|, |v_{q}(a)|\le p^{\,r-1/2}$ then
\[
\bigl|A_{3}(p^{r})\bigr|\,\le\,\frac{q\sqrt{q}}{\varphi^{3}(q)}\cdot 3\varphi(q)\,=\,\frac{3q\sqrt{q}}{\varphi^{2}(q)}
\,=\,\frac{3p^{\frac{\scriptstyle 3r}{\scriptstyle 2\mathstrut}}}{p^{2(r-1)}(p-1)^{2}}\,=\,
3\Bigl(\frac{p}{p-1}\Bigr)^{2}p^{-\,\frac{\scriptstyle r}{\scriptstyle 2\mathstrut}}
\]
for even $r$ and
\[
\bigl|A_{3}(p^{r})\bigr|\,\le\,\frac{q\sqrt{q}}{\varphi^{3}(q)}\cdot 3p^{r-\frac{\scriptstyle 1}{\scriptstyle 2\mathstrut}}\,=\,
\frac{3p^{\frac{\scriptstyle 3r}{\scriptstyle 2\mathstrut}}\cdot p^{r-\frac{\scriptstyle 1}{\scriptstyle 2\mathstrut}}}{p^{3(r-1)}(p-1)^{3}}\,=\,
3\Bigl(\frac{p}{p-1}\Bigr)^{3}p^{-\,\frac{\scriptstyle r+1}{\scriptstyle 2\mathstrut}}
\]
for odd $r$. Thus, if $n = 2h$ is even then
\begin{multline*}
\biggl|\,\sum\limits_{r=n}^{s}A_{3}(p^{r})\biggr|\,\le\,\sum\limits_{\substack{r\ge 2h \\ r\equiv 0\pmod{2}}}3\Bigl(\frac{p}{p-1}\Bigr)^{2}p^{-\,\frac{\scriptstyle r}{\scriptstyle 2\mathstrut}}\,+\,
\sum\limits_{\substack{r\ge 2h+1 \\ r\equiv 1\pmod{2}}}3\Bigl(\frac{p}{p-1}\Bigr)^{3}p^{-\,\frac{\scriptstyle r+1}{\scriptstyle 2\mathstrut}}\,=\\
=\,3\Bigl(\frac{p}{p-1}\Bigr)^{2}\sum\limits_{\ell = h}^{+\infty}p^{-\ell}\,+\,3\Bigl(\frac{p}{p-1}\Bigr)^{3}
\sum\limits_{\ell = h}^{+\infty}p^{-\ell-1}\,=\\
=\,3\Bigl(\frac{p}{p-1}\Bigr)^{3}p^{-h}\,+\,3\Bigl(\frac{p}{p-1}\Bigr)^{3}\frac{p^{-h}}{p-1}\,=\,3\Bigl(\frac{p}{p-1}\Bigr)^{4}p^{-h}\,=
\,3\Bigl(\frac{p}{p-1}\Bigr)^{4}p^{-\,\frac{\scriptstyle n}{\scriptstyle 2\mathstrut}}.
\end{multline*}
Similarly, if $n = 2h+1$ then
\begin{multline*}
\biggl|\,\sum\limits_{r=n}^{s}A_{3}(p^{r})\biggr|\,\le\,\sum\limits_{\substack{r\ge 2h+1 \\ r\equiv 1\pmod{2}}}3\Bigl(\frac{p}{p-1}\Bigr)^{3}p^{-\,\frac{\scriptstyle r+1}{\scriptstyle 2\mathstrut}}\,+\,
\sum\limits_{\substack{r\ge 2h+2 \\ r\equiv 0\pmod{2}}}3\Bigl(\frac{p}{p-1}\Bigr)^{2}p^{-\,\frac{\scriptstyle r}{\scriptstyle 2\mathstrut}}\,=\\
=\,3\Bigl(\frac{p}{p-1}\Bigr)^{3}\frac{p^{-h}}{p-1}\,+\,3\Bigl(\frac{p}{p-1}\Bigr)^{2}
\frac{p^{-h}}{p-1}\,=\,
\,3\Bigl(2-\frac{1}{p}\Bigr)\Bigl(\frac{p}{p-1}\Bigr)^{4}p^{-\,\frac{\scriptstyle n+1}{\scriptstyle 2\mathstrut}}.
\end{multline*}
Lemma is proved. $\Box$
\vspace{0.3cm}

\textsc{Lemma 6.4.} \emph{Under the conditions of Lemma} 6.3, \emph{one has}
\[
\biggl|\,\sum\limits_{r=n}^{s}A_{4}(p^{r})\biggr|\,\le\,6\Bigl(\frac{p}{p-1}\Bigr)^{4}p^{-\,n}.
\]
\textsc{Proof.} Using the same arguments as above, we see that the formula (\ref{7-2}) for $A_{4}(q)$ contains at most one non-zero term. This leads us to the estimate
\[
|A_{4}(q)|\,\le\,\frac{6q^{2}\varphi(q)}{\varphi^{4}(q)}\,=\,\frac{6q^{2}}{\varphi^{3}(q)}\,=\,\frac{6p^{2r}}{p^{3(r-1)}(p-1)^{3}}\,=\,6\Bigl(\frac{p}{p-1}\Bigr)^{3}p^{-r}.
\]
Summing the last inequality over $n\le r\le s$, we get the desired estimate. $\Box$
\vspace{0.3cm}

\textsc{Lemma 6.5.} \emph{Suppose that $p\ge 3$ is prime and $n\ge 2$, $n\equiv \gamma\pmod{2}$, $\gamma = 0,1$. Then, for any $s\ge n$, the following inequality holds:}
\[
\biggl|\sum\limits_{r=n}^{s}A_{5}(p^{r})\biggr|\,\le\,10\Bigl(\frac{p}{p-1}\Bigr)^{6}\cdot\frac{p^{2}-(-1)^{\gamma}(p-1)}{p^{2\mathstrut}+p+1}\,p^{-\frac{\scriptstyle 3n}{\scriptstyle 2\mathstrut}-\frac{\scriptstyle \gamma}{\scriptstyle 2\mathstrut}}.
\]
\textsc{Proof.} First we proof the estimate
\begin{equation*}
|A_{5}(p^{r})|\,\le\,10\Bigl(\frac{p}{p-1}\Bigr)^{4+\delta}p^{-\,\frac{\scriptstyle 3r}{\scriptstyle 2\mathstrut }-\frac{\scriptstyle \delta}{\scriptstyle 2\mathstrut}},\quad
\delta\,=\,
\begin{cases}
& 0,\quad \text{when}\;\;r\equiv 0\pmod{2},\\
& 1,\quad \text{when}\;\;r\equiv 1\pmod{2}.
\end{cases}
\end{equation*}
Indeed, let $r$ be even. Suppose that there are two non-zero terms in the expressions for $A_{5}(p^{r})$ given by Lemma 5.5
Then the difference of the arguments of corresponding Ramanujan sums $c_{q}$ is divisible  by $p^{r-1}$. But all these differences have the form
\[
\pm\,2^{2}\cdot 5\nu,\quad \pm\,2^{4}\nu,\quad \pm\,2^{2}\cdot 3\nu,\quad \pm\,2^{3}\nu,\quad \pm 2^{2}\nu.
\]
Obviously, this is impossible when $p\ge 7$, $r\ge 2$ or when $p = 3,5$, $r\ge 3$. Hence, in these cases we have
\begin{equation}\label{8-5}
|A_{5}(p^{r})|\,\le\,\frac{q^{2}\sqrt{q}}{\varphi^{5}(q)}\cdot 10\varphi(q)\,=\,\frac{10q^{2}\sqrt{q}}{\varphi^{4}(q)}\,=\,10\Bigl(\frac{p}{p-1}\Bigr)^{4}p^{-\,\frac{\scriptstyle 3r}{\scriptstyle 2\mathstrut}}.
\end{equation}
If $p =5$ and $r = 2$ then $5^{r-1} = 5$ may divide the numbers $m\pm 10\nu$. The corresponding sums $c_{q}$ occur in the expression for $A_{5}(p^{r})$ with coefficients equal to 1. Since $1+1<10$, the estimate (\ref{8-5}) holds true in this case.

If $p = 3$, $r=2$ then $p^{r-1}=3$ can divide only the differences
\[
\pm\bigl\{(m+6\nu)-(m-6\nu)\bigr\},\quad \pm\bigl\{(m+10\nu)-(m-2\nu)\bigr\},\quad \pm\bigl\{(m-10\nu)-(m+2\nu)\bigr\}.
\]
Hence, nonzero terms in the expression of Lemma 5.5 should coincide with one of the sums
\[
5c_{9}(m+6\nu)+5c_{9}(m-6\nu),\quad c_{9}(m+10\nu)+10c_{9}(m-2\nu),\quad c_{9}(m-10\nu)+10c_{9}(m+2\nu).
\]
Obviously, the first sum does not exceed $10\varphi(9)$ in modulus, and therefore this case is covered by (\ref{8-5}). In two other cases
we note that the numbers $m+10\nu$, $m-2\nu$ ($m-10\nu$, $m+2\nu$, respectively) are not divided by $9$ simultaneously: in the opposite case, $\nu$ should be divided by $3$. Hence, the sums $c_{9}(m\pm 10\nu)$, $c_{9}(m\mp 2\nu)$ are not equal to $\varphi(9)$ simultaneously. Hence, there is at least one sum which is equal to $0$ or $(-3)$. Thus, one can check that in all such cases the expression $c_{9}(m\pm 10\nu)+10c_{9}(m\mp 2\nu)$ does not exceed $10\varphi(9)-3<10\varphi(9)$ in modulus. Hence, the inequality (\ref{8-5}) holds true.

Suppose now that $r\ge 3$ is odd. Then, using the same arguments as above, we conclude that the expression for $A_{5}(q)$ contains at most one nonzero term. Therefore,
\begin{equation}\label{8-6}
|A_{5}(p^{r})|\,<\,\frac{q^{2}\sqrt{q}}{\varphi^{5}(q)}\cdot 10p^{r-\frac{\scriptstyle 1}{\scriptstyle 2\mathstrut}}\,=\,\frac{10p^{\frac{\scriptstyle 7r}{\scriptstyle 2\mathstrut}-\frac{\scriptstyle 1}{\scriptstyle 2\mathstrut}}}{p^{5r-5}(p-1)^{5}}\,=\,10\Bigl(\frac{p}{p-1}\Bigr)^{5}p^{-\frac{\scriptstyle 3r}{\scriptstyle 2\mathstrut}-\frac{\scriptstyle 1}{\scriptstyle 2\mathstrut}}
\end{equation}
for any $p$ and $k$ under considering.

Next, for even $n = 2h$ and any $s\ge n$ we have
\begin{multline*}
\biggl|\sum\limits_{r=n}^{s}A_{5}(p^{r})\biggr|\,<\,\sum\limits_{\ell\ge h}10\Bigl(\frac{p}{p-1}\Bigr)^{4}p^{-\frac{\scriptstyle 3}{\scriptstyle 2\mathstrut}\cdot 2\ell}\,+\,
\sum\limits_{\ell\ge h}10\Bigl(\frac{p}{p-1}\Bigr)^{5}p^{-\frac{\scriptstyle 3}{\scriptstyle 2\mathstrut}(2\ell+1)-\frac{\scriptstyle 1}{\scriptstyle 2\mathstrut}}\,=\\
=\,10\Bigl(\frac{p}{p-1}\Bigr)^{4}\sum\limits_{\ell = h}^{+\infty}p^{-3\ell}\,+\,10\Bigl(\frac{p}{p-1}\Bigr)^{5}
\sum\limits_{\ell = h}^{+\infty}p^{-3\ell-2}\,=\\
=\,10\Bigl(\frac{p}{p-1}\Bigr)^{4}\frac{p^{-3h}}{1-p^{-3}}\,+\,10\Bigl(\frac{p}{p-1}\Bigr)^{5}
\frac{p^{-3h-2}}{1-p^{-3}}\,=\,10\Bigl(\frac{p}{p-1}\Bigr)^{6}\,\frac{p^{2}-p+1}{p^{2}+p+1}\,p^{-\frac{\scriptstyle 3n}{\scriptstyle 2\mathstrut}}.
\end{multline*}
Finally, for odd $n = 2h+1$ we get
\begin{multline*}
\biggl|\sum\limits_{r=n}^{s}A_{5}(p^{r})\biggr|\,<\,\sum\limits_{\ell\ge h}10\Bigl(\frac{p}{p-1}\Bigr)^{5}p^{-\frac{\scriptstyle 3}{\scriptstyle 2\mathstrut}(2\ell+1)-\frac{\scriptstyle 1}{\scriptstyle 2\mathstrut}}\,+\,
\sum\limits_{\ell\ge h+1}10\Bigl(\frac{p}{p-1}\Bigr)^{4}p^{-\frac{\scriptstyle 3}{\scriptstyle 2\mathstrut}\cdot 2\ell}\,=\\
=\,10\Bigl(\frac{p}{p-1}\Bigr)^{5}
\sum\limits_{\ell = h}^{+\infty}p^{-3\ell-2}\,+\,10\Bigl(\frac{p}{p-1}\Bigr)^{4}\sum\limits_{\ell = h+1}^{+\infty}p^{-3\ell}\,=\\
=\,\Bigl(\frac{p}{p-1}\Bigr)^{5}\frac{p^{-3h+1}}{p^{3}-1}\,+\,10\Bigl(\frac{p}{p-1}\Bigr)^{4}\frac{p^{-3h}}{p^{3}-1}\,=\,
10\Bigl(\frac{p}{p-1}\Bigr)^{6}\cdot\frac{p^{2}+p-1}{p^{2}+p+1}\,p^{-3h-2}.
\end{multline*}
Since $3h+2 = \displaystyle \frac{3n}{2\mathstrut}+\frac{1}{2\mathstrut}$, the the desired estimate follows. Lemma is proved. $\Box$
\vspace{0.5cm}

\textbf{7. Non-vanishing of singular series $\boldsymbol{\varkappa_{k}(q;a,b,m)}$ for $\boldsymbol{q}$ coprime to $\boldsymbol{6}$.}
\vspace{0.5cm}

Here we state the following assertion.
\vspace{0.3cm}

\textsc{Theorem 7.1.} \emph{Suppose that $q$ is coprime to $6$ and let $k\ge 3$ be any fixed integer. Then, for any triple $(a,b,m)$ with the conditions $1\le a,b,m\le q$, $(ab,q)=1$ the following inequalities hold:}
\begin{align*}
\varkappa_{k}(q;a,b,m)\,\ge\,
\begin{cases}
& C_{1}\exp{\Bigl(-\,\displaystyle \frac{C_{2}\sqrt{\ln{q\mathstrut}}}{\ln\ln{q}}\Bigr)}, \quad \textit{if }\;k=3,\\
& C_{3}(\ln\ln{q})^{-6},\quad \textit{if }\; k=4,\\
& 10^{5}, \quad \textit{for any }\; k\ge 5
\end{cases}
\end{align*}
\emph{where the constants $C_{j}$, $j = 1,2,3$, are absolute.}
\vspace{0.3cm}

\textsc{Corollary.} \emph{Under the conditions of Theorem} 7.1, \emph{the formula for $I_{k}(N)$ of Theorem} 3.1 \emph{is asymptotic.}
\vspace{0.3cm}

To prove this theorem, we need some auxilliary lemmas.
\vspace{0.3cm}

\textsc{Lemma 7.1.} \emph{Suppose that $k\ge 3$ is fixed and let $p\ge 7$ be a prime. Then the inequality}
\[
\varkappa_{k}(q)\,=\,\varkappa_{k}(a,b,m;q)\,>\,c_{1}\,=\,\frac{1}{23}
\]
\emph{holds for any $q = p^{n}$, $n\ge 1$, and for any triple $(a,b,m)$ such that} $1\le a,b,m\le q$, $(ab,p)=1$.
\vspace{0.1cm}

\textsc{Proof.} Using the estimates of Lemmas 6.1 and 6.2 together with the fact that
\begin{equation}\label{9-1}
A_{k}(p^{n}) = 0\quad\text{for}\quad n\ge 2\quad \text{when}\quad e=\biggl(\frac{c}{p}\biggr)\,=\,-1,\quad c\equiv ab\pmod{p},
\end{equation}
we easily get
\begin{multline}\label{9-2}
\bigl|\varkappa_{k}(q)\,-\,1\bigr|\,\le\,|A_{k}(p)|\,+\,\sum\limits_{\nu = 2}^{n}|A_{k}(p^{\nu})|\,\le\\
\le\,\frac{(2\sqrt{p})^{k}}{(p-1)^{k-1}}\,\biggl(\frac{1}{4}+\frac{1}{4p}\biggr)\,+\,(1+e)\biggl(\frac{2p}{p-1}\biggr)^{k-1}\,\frac{p^{2-k}}{1-p^{1-k/2}}\,\le\\
\le\,\frac{(2\sqrt{p})^{k}}{(p-1)^{k-1}}\biggl(\frac{1}{4}+\frac{1}{4p}\,+\,\frac{1}{p^{\,k/2-1}-1}\biggr).
\end{multline}
Denote the right-hand side of (\ref{9-2}) by $h(p;k)$. Then
\begin{align*}
\frac{\partial h}{\partial p}\,=& \,-\,\frac{(2\sqrt{p})^{k}}{(p-1)^{k}}\biggl\{\biggl(\frac{k}{2}\Bigl(1+\frac{1}{p}\Bigr)-1\biggr)
\biggl(\frac{1}{4}+\frac{1}{4p}+\frac{1}{p^{\,k/2-1}-1}\biggr)\,+\\
& +\,\frac{p-1}{p^{2}}\biggl(\frac{1}{4}+\biggl(\frac{k}{2}-1\biggr)\,\frac{p^{\,k/2}}{(p^{\,k/2-1}-1)^{2\mathstrut}}\biggr)\biggr\},\\
\frac{\partial h}{\partial k}\,=& \,-\,\frac{(2\sqrt{p})^{k}}{(p-1)^{k-1}}\biggl\{\ln{\biggl(\frac{p-1}{2\sqrt{p}}\biggr)}\biggl(\frac{1}{4}+\frac{1}{4p}\,+\,\frac{1}{p^{\,k/2-1}-1}\biggr)\,+\,
\frac{1}{2}(\ln{p})\,\frac{p^{\,k/2-1}}{(p^{\,k/2-1}-1)^{2}}\biggr\}.
\end{align*}
One can check that
\[
\frac{\partial h}{\partial p}\,<\,0,\quad \frac{\partial h}{\partial k}\,<\,0
\]
for $p\ge 7$, $k\ge 6$. Hence, $h(p,k)$ decreases in this domain with respect to each variable. The direct calculation shows that $h(7,6)=0.865398\ldots$. Hence,
\[
h(p,k)\,\le\,h(7,6),\quad \varkappa_{k}(q)\,>\,1-h(7,6)\,=\,0.134602\ldots>c_{1}
\]
for any $p\ge 7$, $k\ge 6$ and $q = p^{n}$, $n \ge 1$.

Further, since $\partial h/\partial p<0$ for any $p,k\ge 3$ then the functions $h(p,k)$, $k = 3,4,5$, decreases with respect to $p$. Thus, the equalities $h(11,4)=0.721600\ldots$, $h(23,3)=0.955938\ldots$ imply that
\begin{align*}
& \varkappa_{k}(q)\,>\,1-h(11,4)\,>\,0.2784\ldots>c_{1}\quad \text{for any}\;\; k = 4,5,\;\;p\ge 11,\\
& \varkappa_{k}(q)\,>\,1-h(23,3)\,>\,0.04406\ldots\,>\,c_{1}\quad \text{for any}\;\; k = 3,4,5,\;\; p\ge 23.
\end{align*}
Therefore, it is sufficient to check the assertion of the lemma for (a) $k = 3$, $11\le p\le 19$ and for (b) $p = 7$, $3\le k\le 5$.

Since $S(a,b;p) = S(1,c;p)$ for $c\equiv ab\pmod{p}$ then we have $S(fa,fb;p) = S(1,cf^{2};p)$ for any $f$, $(f,p)=1$. If $f$ runs through the reduced residual system modulo $p$ then $f^{2}c$ runs through all quadratic residues or non-residues depending on whether $c$ is a quadratic residue modulo $p$ or no. In view of (\ref{9-1}), Lemmas 6.3, 6.4 imply the inequality
\begin{equation}\label{9-3}
|\varkappa_{k}(q)-1|\,<\,\frac{2}{(p-1)^{k}}\sum\limits_{\nu}|S(1,\nu;p)|^{k}\,+\,\frac{1}{2}(1+e)\,\Phi_{k}(p),
\end{equation}
where $\nu$ runs through all residues modulo $p$ such that $\Bigl(\displaystyle \frac{\nu}{p}\Bigr)=e$
and
\[
\Phi_{3}(p)\,=\,\frac{3p^{3}}{(p-1)^{4}},\quad \Phi_{4}(p)\,=\,\frac{6p^{2}}{(p-1)^{4}},\quad \Phi_{5}(p)\,=\,\frac{10p^{3}}{(p-1)^{6}}\,\frac{p^{2}-p+1}{p^{2\mathstrut}+p+1}.
\]
The direct calculation shows that
\begin{align*}
&\quad \text{in the case}\;\;e=1: && \quad \text{in the case}\;\;e=-1: \\
& \bigl|\,\varkappa_{4}(7^{n})-1\,\bigr|\,\le\,0.382716\ldots, && \bigl|\,\varkappa_{4}(7^{n})-1\,\bigr|\,\le\,0.868827\ldots,\\
& \bigl|\,\varkappa_{5}(7^{n})-1\,\bigr|\,\le\,0.119817\ldots, && \bigl|\,\varkappa_{5}(7^{n})-1\,\bigr|\,\le\,0.530055\ldots,\\
& \bigl|\,\varkappa_{3}(11^{n})-1\,\bigr|\,\le\,0.860875\ldots, && \bigl|\,\varkappa_{3}(11^{n})-1\,\bigr|\,\le\,0.897571\ldots,\\
& \bigl|\,\varkappa_{3}(13^{n})-1\,\bigr|\,\le\,0.692316\ldots, && \bigl|\,\varkappa_{3}(13^{n})-1\,\bigr|\,\le\,0.812294\ldots,\\
& \bigl|\,\varkappa_{3}(17^{n})-1\,\bigr|\,\le\,0.542971\ldots, && \bigl|\,\varkappa_{3}(17^{n})-1\,\bigr|\,\le\,0.624966\ldots,\\
& \bigl|\,\varkappa_{3}(19^{n})-1\,\bigr|\,\le\,0.499344\ldots, && \bigl|\,\varkappa_{3}(19^{n})-1\,\bigr|\,\le\,0.57458\ldots\;.
\end{align*}
Hence,
\begin{align*}
& \varkappa_{k}(7^{n})\,>\,1-0.868827\ldots \,>\,0.131173\ldots>c_{1}\quad (k = 4,5),\\
& \varkappa_{3}(11^{n})\,>\,1-0.897571\ldots \,>\,0.102428\ldots>c_{1}, \\
& \varkappa_{3}(13^{n})\,>\,1-0.812294\ldots \,>\,0.187706\ldots>c_{1},\\
& \varkappa_{3}(17^{n})\,>\,1-0.624966\ldots \,>\,0.375034\ldots>c_{1},\\
& \varkappa_{3}(19^{n})\,>\,1-0.57458\ldots \,>\,0.42542\ldots>c_{1}\;.
\end{align*}
Therefore, it remains to consider the case $k = 3$, $p = 7$. If $n\ge 2$ then Corollary of Lemma 6.3 implies the bound
\[
\bigl|\varkappa_{3}(7^{\,n})\,-\,\varkappa_{3}(7^{2})\bigr|\,\le\,\biggl|\,\sum\limits_{r = 3}^{n}A_{3}(7^{\,r})\biggr|\,\le\,\frac{3\cdot 7\cdot 13}{6^{4\mathstrut}}\,=\,\frac{91}{432}.
\]
At the same time, the direct calculation based on the formulas of Lemma 5.3 shows that the least value of the quantity
\begin{multline*}
\varkappa_{3}(7^{\,2})\,=\,1\,+\,\frac{1}{6^{3}}\sum\limits_{f=1}^{6}S^{3}(1,(\nu f)^{2},7)\cos{\Bigl(\frac{2\pi mf}{7}\Bigr)}\,+\\
+\,\frac{1}{6^{3}}\bigl(c_{7^{2}}(m+6\nu)\,+\,c_{7^{2}}(m-6\nu)\,+\,3c_{7^{2}}(m+2\nu)\,+\,3c_{7^{2}}(m-2\nu)\bigr)
\end{multline*}
is equal to $\displaystyle \tfrac{7}{9}$  (the minimum is taken over all possible values of $1\le m,\nu\le 7^{\,2}$ where $(\nu,7)=1$; it is attained at $m = 7\mu$, $\mu = 1,\ldots,6$, and any $\nu\not\equiv 0\pmod{7}$). Hence, for any $n\ge 2$ we have
\[
\varkappa_{3}(7^{\,n})\,\ge\,\varkappa_{3}(7^{\,2})\,-\,\frac{91}{432}\,\ge\,\frac{7}{9}\,-\,\frac{91}{432}\,=\,\frac{245}{432} = 0.567129\ldots\;>c_{1}.
\]
Finally, if $n = 1$ then the desired bound follows from the inequalities
\[
\bigl|\varkappa_{3}(7)-1\bigr|\,\le\,\frac{2}{6^{3}}\bigl(|S(1,1;7)|^{3}\,+\,|S(1,2;7)|^{3}\,+\,|S(1,4;7)|^{3}\bigr)\,<\,0.381509\ldots\;.
\]
Lemma is proved. $\Box$
\vspace{0.3cm}

\textsc{Lemma 7.2.} \emph{Suppose that $k\ge 3$ is fixed, $n\ge 1$ and $q = 5^{\,n}$. Then the inequality}
\[
\varkappa_{k}(q)\,=\,\varkappa_{k}(a,b,m;q)\,>\,c_{2}\,=\,\frac{1}{22}
\]
\emph{holds for any triple $(a,b,m)$ with the conditions} $1\le a,b,m\le q$, $(ab,5)=1$.
\vspace{0.1cm}

\textsc{Proof.} Suppose first that $k\ge 4$. Let $\displaystyle e = \Bigl(\frac{ab}{5}\Bigr)$. If $e=-1$ then Lemma 5.1 implies the inequality
\[
\bigl|\varkappa_{k}(q)-1\bigr|\,\le\,|A_{k}(5)|.
\]
If $e=1$ then, by Lemmas 6.4 and 6.1 we have
\[
\bigl|\varkappa_{k}(q)-1\bigr|\,\le\,\bigl|A_{k}(5)\bigr|\,+\,\biggl|\sum\limits_{r=2}^{n}A_{k}(5^{r})\biggr|\,\le\,
\bigl|A_{k}(5)\bigr|\,+\,f(k),
\]
where
\begin{equation*}
f(k)\,=\,
\begin{cases}
\displaystyle \frac{75}{128}, & \text{when}\;\; k = 4,\\
\displaystyle \frac{20\cdot 2^{-\,k}\mathstrut\mathstrut}{1-5^{1-k/2}}, & \text{when}\;\; k\ge 5.
\end{cases}
\end{equation*}
Using the same arguments as in previous lemma, we conclude that
\[
\bigl|A_{k}(5)\bigr|\,\le\,\frac{2}{4^{k}}\sum\limits_{\nu\,:\,(\frac{\scriptstyle \nu}{5}) = e}\bigl|S(1,\nu;5)\bigr|^{k}\,=\,2\bigl(|\xi|^{k}\,+\,|\eta|^{k}\bigr),
\]
where
\[
\xi\,=\,\frac{1}{4}\,S(1,1;5)\,=\,\frac{3-\sqrt{5}}{8},\quad \eta\,=\,\frac{1}{4}\,S(1,4;5)\,=\,\frac{3+\sqrt{5}}{8}
\]
in the case $e = 1$ and
\[
\xi\,=\,\frac{1}{4}\,S(1,2;5)\,=\,-\,\frac{\sqrt{5}+1}{4},\quad \eta\,=\,\frac{1}{4}\,S(1,3;5)\,=\,\frac{\sqrt{5}-1}{4}
\]
in the case $e = -1$. Therefore,
\begin{equation*}
\bigl|\varkappa_{k}(q)-1\bigr|\,<\,h(k;a,b),\quad\text{where}\quad
h(k;a,b)\,=\,
\begin{cases}
\displaystyle & 2\bigl(|\xi|^{k}\,+\,|\eta|^{k}\bigr)\;\; \text{for}\;\; e = -1,\\
\displaystyle & 2\bigl(|\xi|^{k}\,+\,|\eta|^{k}\bigr) + f(k)\;\; \text{for}\;\; e = 1.
\end{cases}
\end{equation*}
If $k\ge 5$ then $h(k;a,b)$ is decreasing function of $k$ in both cases $e = \pm 1$. Hence,
\[
\bigl|\varkappa_{k}(q)-1\bigr|\,<\,\max_{a,b}{h(5;a,b)}\,<\,0.813922\ldots
\]
and therefore $\varkappa_{k}(q)>0.186078\ldots >c_{2}$. Similarly, if $k = 4$ then
\[
\bigl|\varkappa_{4}(q)-1\bigr|\,<\,\max_{a,b}{h(4;a,b)}\,<\,0.953125\ldots,
\]
so we have $\varkappa_{4}(q)>0.046875\ldots > c_{2}$.

Finally, let $k = 3$. Then the direct calculation based on the formulas of Lemma 5.5 shows that the smallest values of the quantities
\begin{align*}
& \varkappa_{3}(5)\,=\,1\,+\,\frac{1}{4^{3}}\sum\limits_{f=1}^{4}S^{3}(1,(\nu f)^{2},5)\cos{\Bigl(\frac{2\pi m f}{5}\Bigr)},\\
& \varkappa_{3}(5^{2})\,=\,\varkappa_{3}(5)\,+\,\frac{1}{4^{3}}\bigl(c_{5^{2}}(m+6\nu)+c_{5^{2}}(m-6\nu)+3c_{5^{2}}(m+2\nu)
+3c_{5^{2}}(m-2\nu)\bigr),\\
& \varkappa_{3}(5^{3})\,=\,\varkappa_{3}(5^{2})\,+\,\frac{1}{4^{3}}\Bigl(\frac{\nu}{5}\Bigr)\bigl(u_{5^{3}}(m+6\nu)+u_{5^{3}}(m-6\nu)+3u_{5^{3}}(m+2\nu)
+3u_{5^{3}}(m-2\nu)\bigr)
\end{align*}
are equal to $\displaystyle \frac{35}{64}, \frac{15}{32}$ and $\displaystyle \frac{15}{32}$, respectively. If $n\ge 4$ then Lemma 6.3 implies the inequality
\[
\bigl|\varkappa_{3}(5^{n})\,-\,\varkappa_{3}(5^{4})\bigr|\,\le\,\frac{27}{256}.
\]
Thus we get
\[
\varkappa_{3}(5^{n})\,\ge\,\frac{15}{32}\,-\,\frac{27}{256}\,=\,\frac{93}{256}.
\]
Lemma is proved. $\Box$
\vspace{0.3cm}

\textsc{Lemma 7.3.} \emph{Let $q$ be any integer coprime to $6$. Then the inequality}
\[
\varkappa_{k}(q)\,=\,\varkappa_{k}(a,b,m;q)\,>\,c
\]
\emph{holds with $c = 10^{-5}$ for any triple} $1\le a,b,m\le q$, $(ab,q)=1$ \emph{and for any} $k\ge 5$.
\vspace{0.3cm}

\textsc{Proof.} Let us split $q$ into the product $q_{1}q_{2}$ where all prime divisors of $q_{1}$ exceed 11 and $q_{2} = 1$ or all prime divisors of $q_{2}$ belong to the set $\{5,7,11\}$.
If $q_{1} = 5^{\alpha}7^{\beta}11^{\gamma}>1$, $\alpha,\beta,\gamma\ge 0$ then Lemmas 7.1 and 7.2 imply
\[
\varkappa_{k}(q_{1})\,=\,\varkappa_{k}(5^{\alpha})\varkappa_{k}(7^{\beta})\varkappa_{k}(11^{\gamma})\,>\,\frac{1}{22}\cdot \frac{1}{23^{2\mathstrut}}.
\]
Further, let $q_{2} = p_{1}^{\alpha_{1}}\ldots p_{s}^{\alpha_{s}}\ne 1$ where $13\le p_{1}<\ldots<p_{s}$. By (\ref{9-2}), for $\nu = 1,\ldots, s$ we have
\[
\varkappa_{k}(p_{\nu}^{\alpha_{\nu}})\,>\,1\,-\,\frac{(2\sqrt{p_{\nu}}\,)^{k}}{(p_{\nu}-1)^{k-1}}\biggl(\frac{1}{4}\Bigl(1+\frac{1}{p_{\nu}}\Bigr)\,+\,\frac{1}{p_{\nu}^{k/2-1}-1}\biggr)\,=\,1-h(p_{\nu},k).
\]
Earlier, we find that the function $h(p,k)$ decreases with respect to both variables $p$ and $k$ in the region $p\ge 13, k\ge 5$. Hence, $0<h(p,k)\le h(13,5) = 0.273667\ldots$
for any such pair. Thus, one can easily check that
\[
\ln{(1-h(p,k))}>-\,\frac{6}{5}h(p,k)
\]
and
\[
h(p,5)\,=\,\frac{32p^{\,5/2}}{p^{4}}\biggl(\frac{p}{p-1}\biggr)^{4}\biggl\{\frac{1}{4}\biggl(1+\frac{1}{p}\biggr)+\frac{1}{p^{3/2\mathstrut}-1}\biggr\}\,<\,\frac{13}{p^{\,3/2}}.
\]
Therefore,
\begin{multline*}
\prod\limits_{p\ge 13}\bigl(1-h(p,k)\bigr)\,\ge\,\prod\limits_{p\ge 13}\bigl(1-h(p,5)\bigr)\,=\,\exp{\Bigl(\sum\limits_{p\ge 13}\ln{(1-h(p,5))}\Bigr)}\,>\\
>\,\exp{\Bigl(\,-\,\frac{6\cdot 13}{5}\sum\limits_{p\ge 13}\frac{1}{p^{\,3/2}}\Bigr)}\,>\,\frac{1}{8}.
\end{multline*}
Finally we get
\[
\varkappa_{k}(q)\,>\,\frac{1}{8\cdot 22\cdot 23^{2\mathstrut}}\,>\,10^{-5}.
\]
Lemma is proved. $\Box$
\vspace{0.3cm}

\textsc{Lemma 7.4.} \emph{There exist absolute constants $C_{1}, C_{2}, C_{3}$ such that the inequalities}
\[
\varkappa_{3}(q)\,>\,C_{1}\exp{\biggl(\,-\,\frac{C_{2}\sqrt{\ln{q}\mathstrut}}{\ln\ln{q}}\biggr)},\quad
\varkappa_{4}(q)\,>\,\frac{C_{3}}{(\ln\ln{q})^{6}}
\]
\emph{holds for any $q$ coprime to $6$ and for any triple} $1\le a,b,m\le q$, $(ab,q)=1$.
\vspace{0.3cm}

\textsc{Proof.} Setting
\[
q_{1} = \prod\limits_{p^{\,\alpha\mathstrut}||q,\;p\le 19}p^{\alpha}, \quad q_{2} = \prod\limits_{p^{\,\alpha\mathstrut}||q,\;p\ge 23}p^{\alpha},
\]
by Lemmas 7.1 and 7.2 we find that
\[
\varkappa_{k}(q_{1})\,>\,22^{-1}\cdot 23^{-5}
\]
for $k = 3,4$. Next, one can check that
\[
h(p,3)\,<\,\frac{4.6}{\sqrt{p}},\quad h(p,4)\,<\,\frac{6}{p}
\]
for any $p\ge 23$. Then
\begin{multline*}
\varkappa_{3}(q_{2})\,\ge\,\exp{\Bigl(\sum\limits_{p|q}\ln{\Bigl(1-\frac{4.6}{\sqrt{p}}\Bigr)}\Bigr)}\,\gg\,\exp{\Bigl(-\,4.6\sum\limits_{p|q}\frac{1}{\sqrt{p}}\Bigr)},\\
\varkappa_{4}(q_{2})\,\ge\,\exp{\Bigl(\sum\limits_{p|q}\ln{\Bigl(1-\frac{6}{p}\Bigr)}\Bigr)}\,\gg\,\exp{\Bigl(-\,6\sum\limits_{p|q}\frac{1}{p}\Bigr)}.
\end{multline*}
Now it remains to note that
\[
\sum\limits_{p|q}\frac{1}{\sqrt{p}}\,\ll\,\frac{\sqrt{\ln{q}\mathstrut}}{\ln\ln{q}},\quad \sum\limits_{p|q}\frac{1}{p}\,\le\,\ln\ln\ln{q}+O(1).
\]
Lemma is proved. $\Box$
\vspace{0.3cm}

Now Theorem 7.1 follows immediately from Lemmas 7.3 and 7.4.
\vspace{0.5cm}

\renewcommand{\refname}{\normalsize{References}}

Maris E.~Changa

Moscow Pedagogical State University,

Institute of Mathematics and Informatics,

107140, Moscow, Krasnoprudnaya str., 14

e\,-mail: \texttt{maris\_changa@mail.ru}

\vspace{0.5cm}

Мaxim А.~Korolev

Steklov Mathematical Institute of Russian Academy of Sciences,

119991, Moscow, Gubkina str., 8

e\,-mail: \texttt{korolevma@mi-ras.ru}

\end{document}